\documentclass[12pt]{amsart}
\usepackage{latexsym, amssymb, amsfonts}
\newenvironment{pf}{\proof[\proofname]}{\endproof}
\theoremstyle{plain}
\newtheorem{Thm}{Theorem}[section]
\newtheorem{Cor}[Thm]{Corollary}

\newtheorem{Mainf}[Thm]{First Main Theorem}
\newtheorem{Mains}[Thm]{Second Main Theorem}

\newtheorem{Prop}[Thm]{Proposition}
\newtheorem{Lem}[Thm]{Lemma}
\newtheorem{Cl}[Thm]{Claim}

\theoremstyle{definition}

\newtheorem{Rem}[Thm]{Remark}

\newtheorem{Emp}[Thm]{}

\numberwithin{equation}{section}

\newcommand{\B}[1]{\Bbb#1}
\newcommand{\cal}[1]{\mathcal{#1}}
\newcommand{\C}[1]{\cal#1}
\newcommand{\pgr}{\text{$\PGU{d-1,1}(\B{R})$}}
\newcommand{\gr}{\text{$\GU{d-1,1}(\B{R})$}}

\newcommand{\ff}[1]{\text {$\cal F(#1)$}}
\newcommand{\isom}{\overset {\thicksim}{\to}}
\newcommand{\oom}[2]{\text{$\Omega^{#1}_{#2}$}}

\newcommand{\om}{\prod_{i=1}^r\oom{d}{K_{w_i}}\otimes_{K_{w_i}}E_v}
\newcommand{\omn}{\prod_{i=1}^r(\oom{d}{K_{w_i}}
\widehat{\otimes}_{K_{w_i}}\widehat{K}_{w_i}^{nr})\otimes_{K_{w_i}}E_v}
\newcommand{\si}[2]{\text{$\Sigma^{#1}_{#2}$}}

\newcommand{\SI}{\prod_{i=1}^r\si{d}{K_{w_i}}\otimes_{K_{w_i}}E_v}
\newcommand{\SIN}{\prod_{i=1}^r\si{d,n}{K_{w_i}}\otimes_{K_{w_i}}E_v}

\newcommand{\surj}{\twoheadrightarrow}

\newcommand{\lla}{\longleftarrow}
\newcommand{\hra}{\hookrightarrow}
\newcommand{\wt}{\widetilde}
\newcommand{\Gm}{\Gamma}
\newcommand{\bo}{\bold}
\newcommand{\bb }{B^{d-1}}
\newcommand{\bbr }{(\bb)^r  }
\newcommand{\g}{\C{G}}
\newcommand{\h}{\frak{h}}
\newcommand{\gm}{\gamma}
\newcommand{\dt}{\delta}
\newcommand{\Dt}{\Delta}
\newcommand{\bs}{\backslash}

\newcommand{\m}{^{\times}}

\newcommand{\ii}{^{int}}
\newcommand{\un}{^{un}}
\newcommand{\op}{^{opp}}
\newcommand{\disc}{^{disc}}
\newcommand{\al}{\alpha}
\newcommand{\af}{\B{A}_{F}^{f}}
\newcommand{\afv}{\B{A}^{f;v_1,...,v_r}_{F}}

\newcommand{\an}{^{an}}
\newcommand{\tw}{^{tw}}
\newcommand{\rs}[1]{Section \ref{S:#1}}
\newcommand{\rl}[1]{Lemma \ref{L:#1}}

\newcommand{\rp}[1]{Proposition \ref{P:#1}}
\newcommand{\rr}[1]{Remark \ref{R:#1}}

\newcommand{\re}[1]{\ref{E:#1}}
\newcommand{\rco}[1]{Corollary \ref{C:#1}}

\newcommand{\sm}{\smallsetminus}
\newcommand{\be}{\infty}

\newcommand{\GU}[1]{\bold{GU_{#1}}}

\newcommand{\PGU}[1]{\bold{PGU_{#1}}}
\newcommand{\GL}[1]{\bold{GL_{#1}}}
\newcommand{\PGL}[1]{\bold{PGL_{#1}}}

\newcommand{\PSL}[1]{\bold{PSL_{#1}}}
\newcommand{\Mat}{\text{Mat}}

\newcommand{\pr}{\text{pr}}
\newcommand{\Ker}{\text{Ker}}
\newcommand{\val}{\text{val}}

\newcommand{\diag}{\text{diag}}
\newcommand{\Spec}{\text{Spec}}
\newcommand{\Aut}{\text{Aut}}

\newcommand{\Inn}{\text{Inn}}
\newcommand{\Ad}{\text{Ad}}

\newcommand{\Gal}{\text{Gal}}

\newcommand{\Fr}{\text{Fr}}

\newcommand{\PG}{\bold{PG}}
\newcommand{\PH}{\bold{PH}}
\newcommand{\Lie}{\text{Lie}}

\newcommand{\Hom}{\text{Hom}}

\begin{document}


\title[$P$-adic uniformization II]%
{$P$-adic uniformization\\
   of Unitary Shimura Varieties II.}
\author[Yakov Varshavsky]{Yakov Varshavsky}
\address{School of Mathematics, Fuld Hall\\ 
Institute for Advanced Study\\ 
Olden Lane, Princeton\\ 
NJ 08540, USA}
\email{vyakov@math.toronto.edu}
\curraddr{Department of Mathematics, University of Toronto, 
100 St. George St., Toronto, Ontario M5S 3G3, CANADA}
\date{\today}
\maketitle
\begin{abstract}
In this paper we show that certain Shimura varieties, uniformized by the product 
of complex unit balls, can be $p$-adically uniformized by the product
(of equivariant coverings) of Drinfeld upper half-spaces
and their equivariant coverings. We also extend a $p$-adic uniformization to automorphic vector bundles.
 It is a continuation of our previous work \cite{V} and contains all cases (up to a central modification) of a 
uniformization by known $p$-adic symmetric spaces. The idea of the proof is to show that an arithmetic
quotient of the product of Drinfeld upper half-spaces cannot be anything else than a certain unitary 
Shimura variety. Moreover, we show that difficult theorems of Yau and Kottwitz appearing in \cite{V} 
may be avoided. 
\end{abstract}

\section{Introduction} \label{S:intro}

Let $M$ be a Hermitian symmetric domain (=Hermitian symmetric space of non-compact type), 
and let $\Dt$ be a  torsion-free cocompact lattice
in $\Aut(M)$. Then the quotient $\Dt\bs M$ is  a complex manifold, which has a 
unique structure of a complex projective variety $Y_{\Dt}$ 
(see \cite[Ch. IX, $\S$3]{Sh}). 
A well-known theorem says that when $\Dt$  is an arithmetic congruence subgroup, 
the Shimura variety $Y_{\Dt}$ has a canonical structure over some number field
$E$ (see for example \cite[II, Thm. 5.5]{Mi1}). 

Let $v$ be a prime of $E$. We are interested in a question whether $Y_{\Dt}$ can be
$p$-adically (or more precisely $v$-adically) uniformized. By this we mean that the $E_v$-analytic space 
$(Y_{\Dt}\otimes_E E_v)\an$ is isomorphic to $\Dt\bs\Omega$ for some $E_v$-analytic
symmetric space $\Omega$ and some arithmetic group $\Dt$, acting on $\Omega$  discretely.
In the cases where a $p$-adic uniformization exists we are interested in the relation between 
$M$ and $\Omega$, $\Dt$ and $\Gm$.

The main obstacle for attacking such a general problem is that there is no general definition
of a $p$-adic symmetric space. 
The only $p$-adic analytic spaces which are generally called ``symmetric''
can be described as follows. For a $p$-adic field $L$ and a natural number $d\geq 2$ 
let $\Omega_L^d$ be the open analytic subspace of $\B{P}^{d-1}_L$ obtained by removing the union of all
$L$-rational hyperplanes.

 $\Omega_L^d$ is called $(d-1)$-dimensional Drinfeld upper half-space over $L$ and 
has a lot of good properties. For example, it is a generic fiber of some explicitly
constructed regular formal scheme
$\widehat{\Omega}_L^d$, flat over the ring of integers $\C{O}_L$ of $L$. Moreover,
the natural $L$-rational action of $\PGL{d}(L)$ on $\Omega_L^d$ extends to an 
$\C{O}_L$-rational action on $\widehat{\Omega}_L^d$
(see \cite{Mus} or \cite{Ku}).
Both $\widehat{\Omega}_L^d$ and $\Omega_L^d$ are closely connected with the Bruhat-Tits
building  of $\bo{SL_d}(L)$ (see \cite{Mus} and \cite{Be}), and   $\PGL{d}(L)$ is the group of all automorphisms of  
$\widehat{\Omega}_L^d$ over $\C{O}_L$
(see \cite[Prop. 4.2]{Mus}) and of $\Omega_L^d$ over $L$ (see \cite{Be}).
Furthermore,  Drinfeld constructed a projective system $\{\Sigma_L^{d,n}\}_{n\in\B{N}}$ of 
$\GL{d}(L)$-equivariant finite \'etale Galois coverings of $\Omega_L^d\widehat{\otimes}_L\widehat{L}^{nr}$,
where  $\widehat{L}^{nr}$ is the completion of the maximal unramified extension of $L$
(see \cite{Dr} or \cite[1.4.1]{V}). 
(For simplicity of notation we introduce a pro-analytic space 
$\Sigma^d_L:=\underset{\underset{n}{\lla}}{\lim}\,\Sigma_L^{d,n}$, which was written as 
$\{\Sigma_L^{d,n}\}_n$ in the notation of \cite[Def. 1.3.3-1.3.5]{V}).

But the most important for us property of $\Omega_L^d$  is that for each torsion-free cocompact lattice $\Gm\subset\PGL{d}(L)$ the quotient
$\Gm\bs\Omega_L^d$ exists and has a unique structure of a smooth projective variety $X_{\Gm}$ over $L$.
Moreover, $X_{\Gm}$ satisfies the Hirzebruch proportionality principle, and its canonical divisor is ample
(see \rl{quot} and \rp{Hirz} for a generalization for products). These and some other properties of the
$\Omega_L^d$'s (and of the $\Sigma_L^d$'s) enable us to call them ``symmetric spaces''.

Since there are essentially no other known examples of $p$-adic analytic spaces having analogous properties, especially
the algebraization property of the quotients (compare \cite[Thm 4.8]{PV}), 
we have to restrict ourself to Shimura varieties, 
$p$-adically uniformized by the products of (the equivariant coverings of) Drinfeld upper
half-spaces. This enables us to obtain some arithmetic information about these
Shimura varieties such as a description of their reduction modulo $p$. (To be more precise the symmetric
spaces of \cite[Thm 4.8]{PV} include also some equivariant projective bundles over Drinfeld upper half-spaces.
However our Second Main Theorem implies that these spaces $p$-adically uniformize the corresponding
 projective bundles over Shimura varieties.)

For the proof of a $p$-adic uniformization there are two completely different methods. The first one, due to Cherednik 
\cite{Ch}, is based on Ihara's method of elliptic elements. Cherednik treated the 
case of a uniformization by the $\Omega^2$'s of Shimura curves associated to 
adjoint groups. 

The second method, due to Drinfeld \cite{Dr} (explained in more detail in \cite{BC}), is based on constructing
moduli problems whose solution are simultaneously  Shimura variety and a $p$-adically uniformized variety.
Drinfeld treated the case of a uniformization by the $\Sigma^2$'s, but only for Shimura curves defined over $\B{Q}$.
Later Rapoport and Zink developed Drinfeld's method (see \cite{RZ1, Ra}),  and their recent work \cite{RZ2}
 treats the case of a uniformization by the products of the $\Omega_L^d$'s 
(and of the $\Sigma_L^d$'s). We also would like to mention the work of Boutot and Zink \cite{BZ}, where they 
deduce the uniformization for curves from \cite{RZ2} (compare Section 5).  

Generalizing Cherednik's method, the author treated in his thesis \cite{V} the case of Shimura varieties,
 uniformized by (the equivariant coverings of) Drinfeld upper half-spaces.
 Moreover, the result was extended to the case of automorphic vector bundles and standard principal bundles.

As was indicated in  the introduction of \cite{V}, our method carries to the case of products. 
We do this here, and we also show that the
use in \cite{V} of the classification of algebraic groups and of difficult theorems of Kottwitz and 
Yau was unnecessary. The main new idea of the proof in comparison to that of \cite{V} is to work not 
just with a Shimura variety but with a triple consisting of a Shimura variety,
a standard principal bundle with a flat connection over it, and an equivariant map from the standard
principal bundle to the corresponding Grassmann variety.

 For the convenience of the reader we now compare the results of this paper with those of 
\cite{RZ2}. Note first that our Shimura varieties differ from those of \cite{RZ2} by some abelian 
twist (see \rr{comp}). The Shimura varieties treated by  \cite{RZ2} are moduli varieties of 
abelian varieties with some additional structure, whereas our choice enabled us to write 
the p-adic uniformization in a simpler form (without the twist appearing in \cite[Prop. 6.49]{RZ2}). 
Moreover, our Shimura varieties have weight defined over $\B{Q}$ (the trivial one), 
therefore they are moduli varieties of abelian motives (see \cite[Thm. 3.3.1]{Mi4}) 
and have some other advantages (see for example \cite[III, \S8]{Mi1}).
 
Our result is almost identical to (the twist of) that of \cite{RZ2}, but a more technically involved 
modular approach of Rapoport and Zink requires them to make some unnecessary assumptions 
(see \cite[6.38]{RZ2}), while our group-theoretic method applies without any difference to 
the general case. We should mention that our proof also makes an indirect use of the moduli 
theory of $p$-divisible groups, since this is the only known approach allowing to construct 
 the $\Sigma_L^d$'s. However, for most of the applications one is only interested in Shimura varieties, 
which are maximal at $p$, and these Shimura varieties are $p$-adically uniformized 
by the maximal unramified extension of scalars of the product of the $\Omega_L^d$'s, 
whose construction is completely elementary.

The paper is organized as follows. In the second and the third sections we formulate and prove our main 
results respectively. As an application we find an algebraic connection between Drinfeld upper half-spaces
and the complex unit balls. In the fourth section we give an alternative more direct approach to the 
differential-geometric part of the proof, which uses Yau's theorem on the existence and the uniqueness
of the K\"ahler-Einstein metric. In the fifth section we deduce the analog of our Main Theorems for 
quaternion Shimura varieties  
from the unitary case. In the last section we show that the Shimura varieties treated in this paper
are the most general ones (up to a central modification) which can be $p$-adically uniformized by the product of 
Drinfeld upper half-spaces. The appendix treats the Hirzebruch proportionality principle for quotients of  
the products of Drinfeld upper half-spaces. It is a generalization of \cite[Thm. 2.2.8]{Ku}.

\centerline{\bf Notation and conventions}

1) For a group $G$  let $Z(G)$ be the center of $G$, and let $PG:=G/Z(G)$
be the adjoint group of $G$.

2) For a Lie group $G$ let $G^0$ be its connected component of the identity.

3) For a totally disconnected topological group $E$ let $\ff{E}$ be the set of all compact and 
open subgroups of $E$, and let $E^{disc}$ be the group $E$ with the discrete
topology.

4) For a subgroup $\Dt$ of a topological group $G$ let $\overline{\Dt}$ be the
closure of $\Dt$ in $G$.

5) For an element $g$ of a group $G$ let $\Inn(g)$ be the inner automorphism of $G$ by $g$.

6) For a natural number $n$ let $I_n$ be the $n\times n$ identity matrix, and let $B^n\subset\B{C}^n$ be 
the $n$-dimensional complex unit ball.
 
7) For an analytic space or a scheme $X$ let $T(X)$ be the tangent bundle on $X$.

8) For a vector bundle $V$ on $X$ and a point $x\in X$ let $V_x$ be the 
fiber of $V$ over $x$.

9) For an algebra $D$ let $D\op$ be the opposite algebra of $D$.

10) For a number field $F$ and a finite set $N$ of finite primes of $F$ let $\af$ be the ring of finite adeles of $F$, 
and let $\B{A}_F^{f;N}$ be the ring of finite adeles of $F$ without the components from $N$.

11) For a field extension $K/F$ let $\bo{R}_{K/F}$ be the functor of the restriction of scalars from $K$ to $F$.

12) For a scheme $X$ over a field $K$ and a field extension $L$ of $K$ we will write $X_L$ or $X\otimes_K L$ instead of
$X\times_{\Spec K}\Spec L$.

13) For a complex analytic space $M$ let $\Aut(M)$ be the topological group of holomorphic automorphisms
equipped with the compact-open topology.
 
14) For an analytic space $X$ over a complete non-archimedean field $K$ and a  complete non-archimedean
field extension $L$ of $K$ let $X\widehat{\otimes}_K L$ be the image of $X$ under the ground field 
extension functor from $K$ to $L$. (The completion sign will be omitted in the case of a finite extension).

15) By a $p$-adic field we mean a finite field extension of $\B{Q}_p$ for some prime number $p$.

16) By a $p$-adic analytic space we mean an analytic space over a $p$-adic field in the sense of Berkovich.
 
17) For a vareity $X$ over afield of complex numbers or over a $p$-adic field let $X\an$ the the 
corresponding complex or $p$-adic analytic space.

\centerline{\bf Acknowledgements}
 
First of all the author wants to thank R. Livn\'e for formulating the problem, for 
suggesting the method of the proof and for his attention and help during all 
stages of the work. He also made his corrections of an earlier version of this 
paper.

The author is also grateful to V. Berkovich for his help on $p$-adic analytic
spaces, to J. Rogawski for the reference to \cite{Cl}, to M. Harris for his suggestion to generalize 
the result to its present form, to A. Nair for the reference to \cite{NR}
and to the referee for useful suggestions.

The work on the paper has begun while the author was a Ph.D. student in the Hebrew University of Jerusalem
and continued iat the Institute for the Advanced Study while the author was supported by the 
NSF grant DMS 9304580.
 
\section{Statement of the results} \label{S:stat}

\begin{Emp} \label{E:not1}

Let $F$ be a totally real number field of degree $g$ over $\B{Q}$, 
let $K$ be a totally imaginary quadratic extension $F$, and let 
$\be_1,...,\be_g$ be the archimedean completions of $F$.
Let $d\geq 2$ and $1\leq r\leq g$ be natural numbers. Let $D\ii$ 
be a central simple algebra over $K$ of dimension $d^2$ 
with an involution of the second kind $\al\ii$ over $F$.
We will make the following assumptions:

i) $\al\ii$ has signature $(d-1,1)$ at $\be_1,...,\be_r$; 

ii) $\al\ii$ is positive definite at $\be_{r+1},...,\be_g$. 

\noindent Let 
$\bo{G\ii}:=\bo{GU}(D\ii, \al\ii)$ be the algebraic group over $F$ of unitary similitudes characterized by
$\bo{G\ii}(R)=\{d\in (D\ii\otimes_F R)\m\,|\,d\cdot\al\ii(d)\in R\m\}$ for each $F$-algebra $R$.

\end{Emp}

\begin{Emp} \label{E:pair}

 Set $\bo{H\ii}:=\bo{R}_{F/\B{Q}}\bo{G\ii}$. For each $i=1,...,g$ fix an embedding $\wt{\be}_i:K\hra\B{C}$, 
extending
$\be_i:F\hra\B{R}$. Then $\wt{\be}_1,...,\wt{\be}_g$ define an identification of
$\bo{H\ii}(\B{R})$ with $\GU{d-1,1}(\B{R})^r\times\GU{d}(\B{R})^{g-r}\subset\GL{d}(\B{C})^g$, which is unique
up to an inner automorphism.

Define a homomorphism $h:\bo{R}_{\B{C}/\B{R}}\B{G}_m\to\bo{H\ii}_{\B{R}}$ by requiring for each $z\in\B{C}\m$
$$h(z)=(\underbrace{\diag(1,...,1,z/\bar{z})^{-1};...;\diag(1,...,1,z/\bar{z})^{-1}}_{r};\underbrace{I_d;...;I_d}_{g-r})
\in\GU{d-1,1}(\B{R})^r\times\GU{d}(\B{R})^{g-r}$$
(compare \rr{comp}). Let $M\ii$ be the $\bo{H\ii}(\B{R})$-conjugacy class of $h$. Then $M\ii$ is isomorphic to
$(B^{d-1})^r:=\underbrace{\bb\times...\times\bb}_{r}$ if $d>2$ and to 
$\h^r:=\underbrace{\h\times...\times\h}_{r}$,
where $\h=\B{C}\sm\B{R}$, if $d=2$. In particular, each connected component of $M\ii$ is isomorphic to $\bbr$. Notice that the pair $(\bo{H\ii}, M\ii)$ satisfies Deligne's axioms (see \cite[1.5  and 2.1]{De} or \cite[II, 2.1]{Mi1}).

\end{Emp}

\begin{Emp} \label{E:Shim}

Let $\wt{X}\ii$ be the canonical model of the Shimura variety corresponding to $(\bo{H\ii}, M\ii)$. 
(For the definition of the canonical model we take that of \cite{Mi3}, which has a different sign convention from those of
\cite{De2} and \cite{Mi1} (see the discussion in \cite[1.10]{Mi3}, \cite[Rem. 3.1.13]{V} and Step 4 of Section 3).
Then $\wt{X}\ii$ is a scheme over the reflex field  
$E\subset\B{C}$ of $(\bo{H\ii}, M\ii)$. Moreover, the group $\bo{G\ii}(\af)$ acts on 
$\wt{X}\ii$ in such a way that for each $S\in\ff{\bo{G\ii}(\af)}$ the quotient $S\bs\wt{X}\ii$ is a projective scheme over $E$ and 
$\wt{X}\ii\cong\underset{\underset{S}{\lla}}{\lim}\,S\bs\wt{X}\ii$.

Set $\bo{PG\ii}(F)_+:=\bo{PG\ii}(F)\cap\bo{PG\ii}(F\otimes_{\B{Q}}\B{R})^0$ and write $\pgr_+$ instead of $\pgr^0$.
Then $(\wt{X}\ii_{\B{C}})\an$, defined as $\underset{\underset{S}{\lla}}{\lim}\,(S\bs\wt{X}_{\B{C}}\ii)\an$,
is isomorphic to 
$$[\bbr\times(\bo{G\ii}(\af)/\overline{\bo{Z(G\ii)}(F)})\disc]/
\bo{PG\ii}(F)_+,$$ 
where $(x,g)\gm=(\gm^{-1}x,g\gm)$ for all $x\in\bbr$, $g\in\bo{G\ii}(\af)$ 
and $\gm\in\bo{PG\ii}(F)_+$. Moreover, the action of $\bo{G\ii}(\af)$ on $\wt{X}\ii$ 
corresponds by this isomorphism to left multiplication on the second factor.

Observe that for each $S\in\ff{\bo{G\ii}(\af)}$
the analytic space $(S\bs\wt{X}\ii_{\B{C}})\an$ is a finite disjoint union of quotients of the form 
$\Dt_{aSa^{-1}}\bs\bbr$ for some $a\in\bo{G\ii}(\af)$, where by 
$\Dt_{aSa^{-1}}\subset\Aut(\bbr)^0\cong\pgr^r_+$ we denote the projection of 
$\bo{G\ii}(F)_+\cap(aSa^{-1})$.
Then each $\Dt_{aSa^{-1}}$ is a cocompact lattice
(torsion-free if $S$ is sufficiently small), hence each geometrically connected component of $S\bs\wt{X}\ii$
is of the form described in the introduction.

\end{Emp}

\begin{Emp} \label{E:refl}

The number field $E\subset\B{C}$ can be described as follows (compare \cite[3.1.1 and Prop. 3.1.3]{V}).
Let $K_0\subset\B{C}$ be the composite of the fields $\wt{\be}_1(K),...,\wt{\be}_r(K)$. Set 
$$\Sigma:=\{\sigma\in\Aut(K_0/\B{Q})\,|\,\forall i=1,..,r \,\exists 
\sigma(i)\in\{1,..,r\}:\sigma(\wt{\be}_i(k))=\wt{\be}_{\sigma(i)}(k)\,\forall k\in K\}.$$ 
Then $E$ is the  subfield of $K_0$, fixed elementwise by $\Sigma$. In particular, the extension
$K_0/E$ is Galois with a Galois group $\Sigma$.

\end{Emp}

\begin{Emp} \label{E:p-ad}

Let $v$ be a finite prime of $E$. Set $X\ii:=\wt{X}\ii\otimes_E E_v$. In this paper we 
are going to show that under certain assumptions $X\ii$ admits a $p$-adic uniformization.

Let $p$ be the restriction of $v$ to $\B{Q}$. Then the completion of 
the algebraic closure of $E_v$ is $\B{C}_p$. Choose a field isomorphism 
$\B{C}\isom\B{C}_p$, extending the natural embedding $E\hra E_v\hra\B{C}_p$.
From now on we identify $\B{C}$ with $\B{C}_p$ by means of this isomorphism.
In particular, we will view $E_v$ as a subfield of $\B{C}$.
 
For each $i=1,...,r$ the embedding $\wt{\be}_i:K\hra\B{C}$ corresponds to an embedding
$\al_i:K\hra\B{C}_p$, which extends to a continuous  embedding
$\al_i:K_{w_i}\hra\B{C}_p$ for some prime $w_i$ of $K$, lying over $p$.
Since the group $\Gal(K_0/E)$ preserves the set $\{\wt{\be}_1,...,\wt{\be}_r\}$,
the set $\{w_1,...,w_r\}$ does not depend on the isomorphism $\B{C}\isom\B{C}_p$.
Let $v_i$ be the restriction of $w_i$ to $F$. To the assumptions made in \re{not1} we add the following:

iii) the $v_i$'s are distinct;

iv)  the  $v_i$'s split in $K$;

\noindent (Notice that these two conditions are satisfied automatically if $p$ splits completely in $K_0$)

v) the Brauer invariant of each $D\ii\otimes_{K} K_{w_i}$ is $1/d$.

\end{Emp}

\begin{Lem} \label{L:comp}
$E_v\subset\B{C}_p$ is the composite field of the $\al_i(K_{w_i})$'s, $i=1,...,r$.
\end{Lem}

\begin{pf}
By the definitions, the closure of the image $K_0\hra\B{C}\isom\B{C}_p$
coincides with the composite field of all the $\al_i(K_{w_i})$'s. Therefore we have to show
that $v$ splits completely in $K_0$. For this it is sufficient to show
that if $\sigma\in\Sigma=\Gal(K_0/E)$ acts continuously on $K_0\subset\B{C}_p$, then $\sigma=1$. 
Such a  $\sigma$ induces a continuous isomorphism $\al_i(K)\isom\al_{\sigma(i)}(K)$ for each $i=1,..,r$.
By the definition this implies that $w_{\sigma(i)}=w_i$. Hence, by our assumption,
${\sigma(i)}=i$, so that $\sigma$ induces the identity on each $\wt{\be}_i(K)$. 
Therefore $\sigma$ acts trivially on the composite field $K_0$. 
\end{pf}

 Now we are going to describe the $p$-adic uniformization of $X\ii$. 

\begin{Emp} \label{E:not2}

Let $D$ be a central simple algebra 
over $K$ of dimension $d^2$ with an involution of the second kind $\al$ over $F$ satisfying the following conditions:

vi) the pairs $(D,\al)$ and $(D\ii, \al\ii)$ are locally isomorphic at all finite places of $F$,
 except at $v_1,...,v_r$;

vii) $D$ splits at $w_1,...,w_r$;

viii) $\al$ is positive definite at all the archimedean places of $F$.

The existence of such $D$ and $\al$ follows from the results of  Kottwitz and Clozel
(see \cite[\S2]{Cl}) as in \cite[Prop. 2.3]{Cl}.

Set $\bo{G}:=\GU{}(D,\al)$, and for each $i=1,...,r$ fix a central skew field $\wt{D}_{w_i}$ over 
$K_{w_i}$ of Brauer invariant $1/d$. 
Set also $\g':=\prod_{i=1}^r F_{v_i}\m\times\bo{G}(\afv)$ and 
$\wt{\g}:=\prod_{i=1}^r \wt{D}_{w_i}\m\times\g'$.
Then for each $i=1,...,r$ the pair consisting of an algebra isomorphism
$D\otimes_K K_{w_i}\isom\Mat_d(K_{w_i})$ (resp. $D\ii\otimes_K K_{w_i}\isom\wt{D}_{w_i}$) 
and the similitude homomorphism $\bo{G}(F_{v_i})\to F_{v_i}\m$ (resp. $\bo{G\ii}(F_{v_i})\to F_{v_i}\m$)
gives us an isomorphism $\bo{G}(F_{v_i})\isom\GL{d}(K_{w_i})\times F_{v_i}\m$ (resp. $\bo{G\ii}(F_{v_i})\isom
\wt{D}_{w_i}\m\times F_{v_i}\m$). Using in addition an algebra isomorphism  
$(D,\al)\otimes_F\afv\isom(D\ii,\al\ii)\otimes_F\afv$ we obtain isomorphisms
$\bo{G}(\af)\isom\prod_{i=1}^r\GL{d}(K_{w_i})\times \g'$ and $\bo{G\ii}(\af)\isom\wt{\g}$.
Abusing notation, we will sometimes write these isomorphisms as equalities.
\end{Emp}

For each $S\in\ff{\g'}$ and each $n\in\B{N}\cup\{0\}$ consider a double quotient
$$\wt{X}_{S,n}=S\bs[\SIN\times\g']/\bo{G}(F)$$
(see \cite[1.3.1 and 1.4.1]{V} for our sign convention, which is different from that of Drinfeld).

\begin{Prop}  \label{P:alg}
Each $\wt{X}_{S,n}$ has a canonical structure of an $E_v$-analytic space and of a projective 
scheme  $X_{S,n}$ over $E_v$.
\end{Prop}

\begin{pf}
The proof will follow closely that of \cite[Prop. 1.5.2]{V}. Set $Z_S=\bo{Z(G)}(F)\cap S$. 
Then the group $Z_S$ is a subgroup of $F\m$, which projects to a subgroup of finite index in 
$\prod_{i=1}^r(E_v\m/\C{O}_{E_v}\m)\cong\B{Z}^r$.
Since each $k\in K_{w_i}\m=Z(\GL{d}(K_{w_i}))$ acts on 
$\Omega^d_{K_{w_i}}\widehat{\otimes}_{K_{w_i}}\widehat{K}_{w_i}^{nr}$ 
by the action of the Frobenius automorphism in the power $d\cdot\val(k)$ on the second factor, 
the quotient $Z_S\bs\omn$ is of the form 
$\prod_{i=1}^r(\Omega^{d}_{K_{w_i}}\otimes_{K_{w_i}}\wt{K}_{w_i})\otimes_{K_{w_i}}E_v$ for some 
finite unramified extensions $\wt{K}_{w_i}$'s of $K_{w_i}$'s. 
Hence this quotient is finite over the $E_v$-analytic space $\om$. 
By the formal arguments of \cite[Prop. 1.5.2]{V} the following lemma completes the proof
\end{pf}

\begin{Lem}  \label{L:quot}
Let $L$ be a $p$-adic field, let $L_1,...,L_r$ be $p$-adic subfields of $L$, let $d\geq 2$ be a
natural number, and let $\Gm$ be a torsion-free cocompact lattice in $\prod_{i=1}^r\PGL{d}(L_i)$. Then

a) the action of $\Gm$ on $\prod_{i=1}^r\Omega^d_{L_i}\otimes_{L_i} L$ is discrete and free;

b) the quotient $\Gm\bs\prod_{i=1}^r\Omega^d_{L_i}\otimes_{L_i} L$ exists and has a unique structure
of a smooth projective variety over $L$, whose canonical bundle is ample.
\end{Lem}

\begin{pf}
 GAGA results imply the uniqueness in b). Rest of the proof follows from the arguments of \cite[Lem. 6]{Be} and \cite[App. 1]{Mus}, working without change in our situation.
\end{pf}

\begin{Emp} \label{E:inv}
Set $X=\underset{\underset{S,n}{\lla}}{\lim}\,X_{S,n}$. Then, as in \cite[Con. 1.5.1]{V}, the group 
$\wt{\g}$ acts on $X$ in such a way that for each $T\in\ff{\wt{\g}}$ the quotient $T\bs X$ is a projective
scheme over $E_v$ and $X=\underset{\underset{T}{\lla}}{\lim}\,T\bs X$.
Moreover,  for each sufficiently small $T\in\ff{\wt{\g}}$ the projection $X\to T\bs X$ is \'etale, and
 the quotient $T\bs X$ is smooth.
\end{Emp}

\begin{Rem} \label{R:cov}
Since the natural projections $(X_{S_1,n_1})\an\to (X_{S_2,n_2})\an$ are \'etale for all $S_1\subset S_2$
sufficiently small and all $n_1\geq n_2$,  we can define, as in \cite[Prop. 1.5.3 f)]{V}, a pro-analytic space 
$X\an:=\underset{\underset{S,n}{\lla}}{\lim}\,X_{S,n}\an$. Then the natural $\wt{\g}$-equivariant 
map $\eta:[\SI\times(\g')\disc]/\bo{G}(F)\to X\an$ of pro-analytic spaces over $E_v$ is \'etale 
and surjective.
\end{Rem}

\begin{Rem} \label{R:unif}
The proof of \rp{alg} (and that of \cite[Prop. 1.5.2]{V}) implies that for each sufficiently small
$T\in\ff{\wt{\g}}$, each connected component of $(T\bs X)\an$ is a finite \'etale cover of an analytic space of
the form $\Gm\bs\prod_{i=1}^r\Omega^d_{K_{w_i}}\otimes_{K_{w_i}} E_v$ for some irreducible arithmetic 
torsion-free cocompact lattice $\Gm\subset\prod_{i=1}^r\PGL{d}(K_{w_i})$. In particular, the canonical bundle
of such an $T\bs X$ is ample.
\end{Rem}

\begin{Mainf}
There exists an isomorphism $\varphi:X\isom X\ii$ commuting with the action of the group
$\wt{\g}=\bo{G\ii}(\af)$.
\end{Mainf}

Now we describe the $p$-adic uniformization of certain automorphic vector bundles and a certain principal bundle.

\begin{Emp} \label{E:bundles}

Let $P\ii$ be the canonical model of the standard $\bo{PH\ii}$-principal
bundle over $X\ii$ (see \cite[Thm. 4.3]{Mi1}). We have $$(P\ii_{\B{C}})\an\cong [\bbr\times(\bo{PH\ii}_{\B{C}})\an
\times(\bo{G\ii}(\af)/\overline{\bo{Z(G\ii)}(F)})\disc]/\bo{PG\ii}(F)_+,$$ 
where $\bo{PG\ii}(F)_+$ acts on $(\bo{PH\ii}_{\B{C}})\an=(\bo{PG\ii}_{F\otimes_{\B{Q}}\B{C}})\an$ 
by right multiplication. Then $P\ii_{\B{C}}$ has a natural $\wt{\g}$-invariant flat connection
$\wt{\C{H}}\ii$  (see \cite[1.9.15-1.9.19]{V} for the definitions) such that the restriction of 
$(\wt{\C{H}}\ii)\an$ to $\bbr\times(\bo{PH\ii}_{\B{C}})\an\times\{1\}$ is trivial, that is 
consists of vectors, tangent to $\bbr$ (compare the proof of \cite[Prop. 4.4.2]{V}). Furthermore,
$P\ii$ is the only $\bo{PH\ii_{E_v}}\times\wt{\g}$-equivariant model of $P\ii_{\B{C}}$,
to which $\wt{\C{H}}\ii$ descends (see for example \cite[4.7]{V}). Denote the descent of $\wt{\C{H}}\ii$ to
 $P\ii$  by $\C{H}\ii$.

Let $\check{M}\ii$ be the Grassmann variety corresponding to the pair $(\bo{H\ii}, M\ii)$ 
(see \cite[III, \S1]{Mi1}). It is a homogeneous space for the group $\bo{PH\ii_E}$,
satisfying $\check{M}\ii_{\B{C}}\cong (\B{P}^{d-1}_{\B{C}})^r$.

Let $\beta_{\B{R}}:\bbr\hra ((\B{P}^{d-1}_{\B{C}})^r)\an$ be the natural (Borel)
embedding, and let $\rho\ii:P\ii_{E_v}\to \check{M}\ii_{E_v}$ be the natural 
$\bo{PH\ii_{E_v}}\times\wt{\g}$-equivariant map (for the trivial action of $\wt{\g}$ on $\check{M}\ii_{E_v}$)
defined by $(\rho\ii_{\B{C}})\an[x,h,g]=h(\beta_{\B{R}}(x))$ for all $x\in\bbr,\;h\in (\bo{PH\ii}_{\B{C}})\an$
and $g\in\wt{\g}$ (see \cite[4.3.2 and 4.3.4]{V}).

Let $W\ii$ be a $\bo{PH\ii_{E_v}}$-equivariant vector bundle on $\check{M}\ii_{E_v}$, and
let $V\ii=V\ii(W\ii)$ be the canonical model of the automorphic vector bundle
on $X\ii$ corresponding to $W\ii$. Then 
$(V\ii_{\B{C}})\an\cong [\beta_{\B{R}}^*(W\ii_{\B{C}})\an\times
(\bo{G\ii}(\af)/\overline{\bo{Z(G\ii)}(F)})\disc]/\bo{PG\ii}(F)_+$ and 
$V\ii\cong \bo{PH\ii_{E_v}}\bs(\rho\ii)^*(W\ii)$
(see \cite[III, Prop. 3.5 and Thm. 5.1]{Mi1}).

\end{Emp}

We now describe the corresponding objects on the $p$-adic side. 

\begin{Emp} \label{E:tw}

Set $\bo{H}:=\bo{R}_{F/\B{Q}}\bo{G}$. Then for some group $\bo{\wt{H}}$ over $E_v$ we have natural
isomorphisms $\bo{PH_{E_v}}\cong(\PGL{d})^r\times\bo{\wt{H}}$ and 
$\bold{PH\ii_{E_v}}\cong\prod_{i=1}^r\bo{PGL_1}(\wt{D}_{w_i})_{\bold{E_v}}
\times\bold{\wt{H}}$, where the first $r$ factors  correspond 
to the natural embeddings $F\hra F_{v_i}\overset{\al_i}{\hra}E_v$. 
Let $\bo{PH_{E_v}}$ acts on 
$(\B{P}^{d-1}_{E_v})^r$ by the natural action of the first $r$ factors 
and the trivial action of the last one.

Let $\pi\in E_v$ be a uniformizer, and let $\wt{\Pi}$ be an element of $\GL{d}(E_v)$ satisfying 
$\wt{\Pi}^d=\pi$. For each $i=1,...,r$ set $d_i:=[E_v:\al_i(K_{w_i})]$.
Finally, denote the  projection of $(\wt{\Pi}^{d_1},...,\wt{\Pi}^{d_r})$ to $\PGL{d}(E_v)^r$ by 
$\Pi'$, and set $\Pi:=(\Pi',1)\in\PGL{d}(E_v)^r\times\bold{\wt{H}}(E_v)\cong\PH(E_v)$.
Let $E_v^{(d)}$ be the unramified field extension of $E_v$ of degree $d$.
Since the Brauer invariant of each $\wt{D}_{w_i}\otimes_{K_{w_i}}E_v$ is $d_i/d$ (see \cite[Ch. VI, Sec. 1,
Thm. 3]{CF}), the group $\bold{PH\ii_{E_v}}$ is 
isomorphic to the quotient of $\bold{PH_{E_v}}\otimes_{E_v}E_v^{(d)}$ by the equivalence
relation $\Fr(x)\sim\Pi^{-1}x\Pi$, where $\Fr\in\Gal(E_v^{(d)}/E_v)$ is the 
Frobenius automorphism.

For an $E_v$-scheme $Y$ with an $E_v$-rational action of the group $\bo{PH_{E_v}}$
 define a twist $Y\tw :=(\Fr(x)\sim\Pi^{-1}x)\bs Y\otimes_{E_v}E_v^{(d)}$. Then 
$Y\otimes_{E_v}E_v^{(d)}\cong Y\tw\otimes_{E_v}E_v^{(d)}$, and the natural action of 
$\bold{PH\ii_{E_v}}$ on $Y\tw$ is $E_v$-rational.
Using the definition of the twist, we see that $\check{M}\ii_{E_v}\cong ((\B{P}^{d-1}_{E_v})^r)\tw$ and 
that there exists a unique
$\bo{PH_{E_v}}$-equivariant vector bundle $W$ on $(\B{P}^{d-1}_{E_v})^r$ such that $W\tw\cong W\ii$.

Let $\beta_v:\SI\surj\om\hra((\B{P}^{d-1}_{E_v})^r)\an$ be the natural map. As in 
\cite[4.1.2-4.1.3 and 4.2.2-4.2.3]{V}, we construct a $\wt{\g}$-equivariant vector bundle
$V=V(W)$ on $X$ and a $\wt{\g}$-equivariant $\bo{PH_{E_v}}$-principal bundle $P$ over $X$
such that $(S\bs V)\an\cong S\bs[\beta_v^*(W)\times\g']/\bo{G}(F)$ and 
$(S\bs P)\an\cong S\bs[\SI\times(\bo{PH_{E_v}})\an\times\g']/\bo{G}(F)$ for each $S\in\ff{\wt{\g}}$.
Of course $V$ and $P$ are unique. $P$ also has (as in the proof of \cite[Prop. 4.4.2]{V})
a natural $\wt{\g}$-invariant flat connection $\C{H}$. As in \cite[Prop. 4.3.2]{V}, there exists a 
$\bo{PH_{E_v}}\times\wt{\g}$-equivariant morphism $\rho:P\to(\B{P}^{d-1}_{E_v})^r$, defined by
$\rho\an[x,h,g]=h(\beta_w(x))$ for each $x\in\SI,\;h\in(\bo{PH_{E_v}})\an$ and $g\in\g'$. 
By the construction, $V$ is naturally isomorphic to the quotient $\bo{PH_{E_v}}\bs\rho^*(W)$. 
\end{Emp}

\begin{Mains} 
a) Any isomorphism $\varphi:X\isom X\ii$ from the First Main Theorem lifts
to a $\wt{\g}$-equivariant isomorphism $\varphi_V:V\isom V\ii$ of automorphic vector bundles.

b) Any $\varphi$ as in a) lifts to a $\wt{\g}$-equivariant isomorphism 
$\varphi_P:P\tw\isom P\ii$ of $\bo{PH\ii_{E_v}}$-principal bundles, 
which maps $\C{H}\tw$ and $\rho\tw$ to $\C{H}\ii$ and $\rho\ii$ respectively.
\end{Mains}

\begin{Rem}
The observations made above show that a) is an immediate consequence of b).
\end{Rem}

The analogs of the Main Theorems for quaternion Shimura varieties will be given in Section 5.

\section{Proof of the Main Theorems} \label{S:proof}


 {\bf Step 1}. Exactly by the same arguments as in  the proof of \cite[Prop. 1.5.3 and 1.6.1]{V} we show that
there exists an inverse limit of $\{(S\bs X_{\B{C}})\an\}_{S\in\ff{\wt{\g}}}$ in the 
category of complex analytic spaces, which we denote by $(X_{\B{C}})\an$. Moreover, the group 
$\wt{\g}$ acts transitively on the set of connected components of $(X_{\B{C}})\an$,
 and $\g_0:=\overline{\bo{Z(G\ii)}(F)}$
is the kernel of the action of the group $\wt{\g}=\bo{G\ii}(\af)$ on $X$. Set $\g:=\wt{\g}/\g_0$. 

Let $M$ be a connected component of  $(X_{\B{C}})\an$, and let
$\Dt$ be the stabilizer of $M$ in $\g$. Then $\Dt$ is naturally embedded into $\Aut(M)\times\g$, and 
$(X_{\B{C}})\an\cong [M\times\g\disc]/\Dt$. Let $\pi:\wt{M}\to M$ be the universal cover of $M$, 
let $\wt{\Dt}\subset\Aut(\wt{M})$ be the set of all automorphisms of $\wt{M}$ which lift some automorphism 
from $\Dt\subset\Aut(M)$, and let $\wt{\pi}:\wt{\Dt}\to\Dt$ be the natural projection. 
Then $(X_{\B{C}})\an\cong [\wt{M}\times\g\disc]/\wt{\Dt}$, where $\wt{\Dt}$ acts on $\g$
though $\wt{\pi}$. 

As in \cite[4.4.1]{V}, the $\g$-equivariant $\bo{PH_{E_v}}$-principal bundle $P$ over 
$X$ with a flat connection $\C{H}$ define a $(\bo{PH}_{\B{C}})\an$-principal bundle $(P_{\B{C}})\an$ over 
$(X_{\B{C}})\an\cong[\wt{M}\times\g\disc]/\wt{\Dt}$ with a flat connection $(\C{H}_{\B{C}})\an$.  
Therefore $\rho:P\to(\B{P}^{d-1}_{E_v})^r$ defines the $\bo{PH}(\B{C})\times\g$-equivariant map 
$(\rho_{\B{C}})\an:(P_{\B{C}})\an\to((\B{P}^{d-1}_{\B{C}})^r)\an$.
Since $\wt{M}$ is simply connected, we show, as in \cite[Prop. 4.4.2]{V}, that there exists a homomorphism 
$j:\wt{\Dt}\to\bo{PH}(\B{C})$ and a  $\bo{PH}(\B{C})\times\g$-equivariant isomorphism  
$(P_{\B{C}})\an\cong [\wt{M}\times(\bo{PH}_{\B{C}})\an\times\g\disc]/\wt{\Dt}$ over $X$,
where $\dt\in\wt{\Dt}$ acts on the second factor by right multiplication by $j(\dt)$. Furthermore, 
 $(\C{H}_{\B{C}})\an$ corresponds then to the natural $\g$-invariant flat connection on the right hand side.

Let $\rho_0$ be the pull-back of $(\rho_{\B{C}})\an$ to 
$\wt{M}\times\{1\}\times\{1\}\cong\wt{M}$, and let 
$i:\wt{\Dt}\to\Aut(\wt{M})\times\bo{PH}(\B{C})\times\g$ be the product of the natural embedding, 
 of $j$ and of $\wt{\pi}:\wt{\Dt}\to\Dt\subset\g$.
It follows that to prove the complex versions of our Main Theorems it will suffice to show 
(see \cite[Lem. 2.2.6 and Rem. 4.4.4]{V}) that 
$M\cong\wt{M}\cong\bbr$, that $\rho_0$ is the Borel embedding, that 
$\Dt\cong\wt{\Dt}\cong\bo{PG\ii}(F)_+$, and that 
$i$ is conjugate to the diagonal embedding of $\bo{PG\ii}(F)_+$ into 
$\prod_{i=1}^r\bo{PG\ii}(F_{\be_i})^0\times\bo{PH\ii}(\B{C})\times(\bo{G\ii}(\af)/\overline{\bo{Z(G\ii)}(F)})$
(for the last property we identify $\wt{M}$ with $\bbr$ and $\wt{\Dt}$ with $\bo{PG\ii}(F)_+$).


 {\bf Step 2}. This differential-geometric part of the proof significantly differs from the 
case of one factor, treated in \cite{V}. The  method we use here is very similar to 
(but was obtained mostly independently of) that of \cite{NR}. Another method will be shown in the next section.

By the definition, $\rho$ induces isomorphisms 
$(\C{H}_{\B{C}_p})_x\isom T_{\rho(x)}(\B{P}^{d-1}_{\B{C}_p})^r$ for all $x\in P$. Therefore
 $\rho_0$ induces isomorphisms between $T_x(\wt{M})$ and $T_{\rho_0(x)}((\B{P}^{d-1}_{\B{C}})^r)\an$ 
for all $x\in\wt{M}$. It follows that $\rho_0$ is a local isomorphism and that
$\rho_0(\wt{M})$ is an open analytic subspace of $((\B{P}^{d-1}_{\B{C}})^r)\an$. 
Furthermore, since $(\rho_{\B{C}})\an$ is $\bo{PH}(\B{C})\times\g$-equivariant, we have 
$\rho_0(\dt x)=j(\dt)\rho_0(x)$ for all $\dt\in\wt{\Dt}$ and $x\in \wt{M}$. In particular,
each $j(\dt)$ induces an automorphism of $\rho_0(\wt{M})$.

Let $j'$ be the composition of $j:\wt{\Dt}\to\PH(\B{C})\cong\PGL{d}(\B{C})^g$ with the projection 
to the product of the first $r$ factors (corresponding to the embeddings $\be_1,...,\be_r$ of $F$ into $\B{R}$). 
Then the local isomorphism $\rho_0$ is $j'$-equivariant.
Denote by $\wt{J}$ the closure of $\wt{\Dt}$ in $\Aut(\wt{M})$.

\begin{Lem} \label{L:extension}
$j'$ can be extended to a continuous homomorphism $\phi:\wt{J}\to\PGL{d}(\B{C})^r$,
and the map $\rho_0$ is $\phi$-equivariant.
\end{Lem}
\begin{pf}

Choose a $g\in\wt{J}$ and a  subsequence $\{\dt_n\}_{n}\subset\wt{\Dt}$, converging to $g$. 
Then for all  $x_1, x_2\in\wt{M}$ the equality $\rho_0(x_1)=\rho_0(x_2)$ implies that
$$\rho_0(g(x_1))=\lim_n\rho_0(\dt_n(x_1))=\lim_n j'(\dt_n)\rho_0(x_1)=\lim_n j'(\dt_n)\rho_0(x_2)=
\lim_n\rho_0(\dt_n(x_2))=\rho_0(g(x_2)).$$
Hence the map $\phi_0(g):\rho_0(x)\mapsto\rho_0(g(x))$ is a well-defined map of $\rho_0(\wt{M})$ 
into itself. Moreover, since $\rho_0$ is a local isomorphism, the constructed map $g\mapsto\phi_0(g)$ is
a continuous homomorphism from $\wt{J}$ to $\Aut(\rho_0(\wt{M}))$, extending $j'$.

Denote the stabilizer of $\rho_0(\wt{M})$ in $\PGL{d}(\B{C})^r$ by $G$. 
Then $G$, being the stabilizer of the closed subset $((\B{P}^{d-1}_{\B{C}})^r)\an\sm\rho_0(\wt{M})$
of $((\B{P}^{d-1}_{\B{C}})^r)\an$, 
is a closed subgroup of $\PGL{d}(\B{C})^r=\Aut((\B{P}^{d-1}_{\B{C}})^r)\an$. 
By the identity theorem (see \cite[Vol. I, A, Thm.3]{Gu}),
the restriction homomorphism $res:G\to\Aut(\rho_0(\wt{M}))$ is injective. Moreover, it induces a
homeomorphism between $G$ and $res(G)\subset\Aut(\rho_0(\wt{M}))$. In particular, the induced topology 
on $res(G)$ is locally compact. Hence $res(G)$ is a closed subgroup of  $\Aut(\rho_0(\wt{M}))$
(see \cite[Prop. 1.4]{Shi}). Therefore $\phi_0(\wt{J})\subset res(G)$, and 
$\phi:=res^{-1}\circ\phi_0:\wt{J}\to G\subset\PGL{d}(\B{C})^r$ is the required homomorphism.
\end{pf}

Our aim is to obtain more information about $\wt{J}$ and $\phi$.
First we show, as in the proof of \cite[Prop. 1.5.3 d), f) and 1.3.8 a)]{V}, that the quotient 
$X':=(\prod_{i=1}^r\wt{D}_{w_i}\m\cdot Z(\wt{\g}))\bs X$ exists, that the projection $X\to X'$ is \'etale,
and that $(X')\an\cong[\prod_{i=1}^r\Omega^d_{K_{w_i}}\otimes_{K_{w_i}} E_v
\times(P\g')^{disc}]/\bo{PG}(F)$. As in \cite[4.2.2-4.2.3]{V}, there exists a unique 
$\bo{PH_{E_v}}$-principal bundle $P'$ over $X'$, satisfying 
$(P')\an\cong[\prod_{i=1}^r\Omega^d_{K_{w_i}}\otimes_{K_{w_i}}E_v\times(\bo{PH_{E_v}})\an\times
(P\g')^{disc}]/\bo{PG}(F)$. Then $P'\cong(\prod_{i=1}^r\wt{D}_{w_i}\m\cdot Z(\wt{\g}))\bs P$, hence
$(P'_{\B{C}})\an\cong[\wt{M}\times(\bo{PH_{\B{C}}})\an\times(P\g')^{disc}]/\wt{\Dt}$
and $(X'_{\B{C}})\an\cong[\wt{M}\times(P\g')^{disc}]/\wt{\Dt}$.
  
For each $x\in\prod_{i=1}^r\Omega^d_{K_{w_i}}(\B{C}_p)$ and $y\in\wt{M}$ set 
$\Gm_x:=\{\gm\in\bo{PG}(F)\,|\,\gm(x)=x\}$ and $\Dt_y:=\{\dt\in\wt{\Dt}\,|\,\dt(y)=y\}$.
Since the projection $X\to X'$ is \'etale, $\wt{M}$ is the universal cover of each connected component
of $(X'_{\B{C}})\an\cong[\wt{M}\times(P\g')^{disc}]/\wt{\Dt}$. Hence each $\Dt_y$ projects injectively to
$P\g'$. In particular, $\Gm_x$'s and $\Dt_y$'s are naturally embedded into 
$\bo{PH}(\B{C}_p)\times P\g'=\bo{PH}(\B{C})\times P\g'$.

The following proposition is a technical heart of the proof (compare \cite[Prop. 2.2.8]{V}).

\begin{Prop} \label{P:ell}
For each $x\in\prod_{i=1}^r\Omega^d_{K_{w_i}}(\B{C}_p)$ there exists a (non-unique) $y\in\wt{M}$ and for 
each  $y\in\wt{M}$ there exists a (non-unique) $x\in\prod_{i=1}^r\Omega^d_{K_{w_i}}(\B{C}_p)$ such that 
the subgroups $\Gm_x$ and $\Dt_y$ are conjugate in $\bo{PH}(\B{C})\times P\g'$. In particular, 
the closure of each 
$j(\Dt_y)\subset\bo{PH}(\B{C})$ is compact, and there exists $y\in\wt{M}$ such that the closure
$\overline{j(\Dt_y)}$ is a maximal compact torus of $\bo{PH}(\B{C})$.
\end{Prop}
\begin{pf}
For each $x\in\prod_{i=1}^r\Omega^d_{K_{w_i}}(\B{C}_p)$ consider a point $[x,1,1]\in
[\prod_{i=1}^r\Omega^d_{K_{w_i}}\otimes_{K_{w_i}}\B{C}_p\times(\bo{PH}_{\B{C}_p})\an\times
(P\g')^{disc}]/\bo{PG}(F)\cong (P'_{\B{C}_p})\an$. It gives us a point of $P'(\B{C}_p)=P'(\B{C})$. 
Therefore it corresponds to a certain point
$[y,h,g]\in[\wt{M}\times(\bo{PH_{\B{C}}})\an\times(P\g')^{disc}]/\wt{\Dt}\cong (P'_{\B{C}})\an$. 
Each $\gm\in\Gm_x\subset\bo{PH}(\B{C})\times P\g'$ fixes $[x,1,1]\in(P'_{\B{C}_p})\an$.
Therefore its conjugate $(h,g)\gm(h,g)^{-1}\in\bo{PH}(\B{C})\times P\g'$ fixes
$[y,1,1]\in[\wt{M}\times(\bo{PH_{\B{C}}})\an\times(P\g')^{disc}]/\wt{\Dt}$. It follows that
 $(h,g)\gm(h,g)^{-1}$ belongs to $\Dt_y$.
The proof of the opposite direction is exactly the same.

By the proven above it will suffice to show the remaining statements for the closures of $\Gm_x$'s
in $\bo{PH}(\B{C})=\bo{PG}(F\otimes_{\B{Q}}\B{C})$. 

As each $\Gm_x$ is contained in $\bo{PH}(\B{R})\cong\PGU{d}(\B{R})^g$, its closure is compact.
Choose an elliptic element $\gm\in\bo{G}(F)\subset D$ with an elliptic point 
$x\in\prod_{i=1}^r\Omega^d_{K_{w_i}}(\B{C}_p)$ (see \cite[Def. 1.7.1]{V}). (Such an element exists,
since by the weak approximation theorem the closure of $\bo{PG}(F)\subset\prod_{i=1}^r\PGL{d}(K_{w_i})$
contains $\prod_{i=1}^r\PSL{d}(K_{w_i})$, and since the set of elliptic elements of 
$\prod_{i=1}^r\PSL{d}(K_{w_i})$ is open and non-empty (see \cite[Prop. 1.7.3]{V}).
By \cite[Lem 1.7.4]{V}, $g$ generates a maximal commutative subfield $L\subset D$, invariant under $\al$. 
Then $\bo{T}:=\bo{H}\cap \bo{R}_{L/\B{Q}}\B{G}_m\subset \bo{R}_{D/\B{Q}}\B{G}_m$ is a maximal 
$\B{Q}$-rational torus of $\bo{H}$, and each $\gm\in\bo{T}(\B{Q})=\bo{G}(F)\cap L\m$ fixes $x$. Moreover,
one check by direct computation that $\bo{T}(\B{Q})$ is the stabilizer of $x$ in $\bo{G}(F)$.
Therefore the projection $\bo{T'}\subset\bo{PH}$ of $\bo{T}$ is also a maximal $\B{Q}$-rational torus
and $\Gm_x=\bo{T'}(\B{Q})$. The theorem on the real approximation (see \cite[(0.4)]{De})
now implies that the closure $\overline{\Gm_x}$ is a maximal compact torus $\bo{T'}(\B{R})$
of $\bo{PH}(\B{C})$.
\end{pf}

Using the proposition, fix a point $y_0\in\wt{M}$ such that the closure of $j'(\Dt_{y_0})$ is 
a maximal compact torus $T$ in $\PGL{d}(\B{C})^r$.

\begin{Emp} \label{E:clos}
Let $J$ be the closure of $j'(\Dt)\subset\PGL{d}(\B{C})^r$. 
Then $\Lie\,J$ is an $\Ad(J)$-invariant real Lie subalgebra of $\Lie\,\PGL{d}(\B{C})^r$. 
Denote by $\Lie\,J_{\B{C}}$ the $\B{C}$-span of 
$\Lie\,J$. As in \cite[Prop. 4.5.2]{V}, the subgroup $j(\Dt)$ is 
Zariski dense in $\bo{PH}_{\B{C}}$, and therefore $\Lie\,J_{\B{C}}$ 
is an ideal in $\Lie\,\PGL{d}(\B{C})^r$. Since  $\Lie\,J_{\B{C}}$ contains the Cartan subalgebra 
$\C{T}:=(\Lie\,T)_{\B{C}}$ of $\Lie\,\PGL{d}(\B{C})^r$, it then has to be all of $\Lie\,\PGL{d}(\B{C})^r$.
Lie algebra $\Lie\,J$ is invariant under the compact group $Ad(T)$, hence it decomposes as
$\prod_{i=1}^{r}\C{J}_i\subset\Lie\,\PGL{d}(\B{C})^r$, and each $\C{J}_i$ is either 
all of $\Lie\,\PGL{d}(\B{C})$ or one of its real form. Thus $J^0$ decomposes as $\prod_{i=1}^{r}J_i$, where
each is $J_i$ either all of $\PGL{d}(\B{C})$ or $\bold{J_i}(\B{R})^0$ for some real form
$\bold{J_i}$ of $\PGL{d}$. In both cases $J^0\cong\bold{J}(\B{R})^0$ for some semi-simple adjoint
real group $\bold{J}$.
\end{Emp}

\begin{Emp} \label{E:surj}
Choose a maximal compact subgroup $J^{comp}\supset T$ of $\PGL{d}(\B{C})^r$. Equip $((\B{P}^{d-1}_{\B{C}})^r)\an$
with a standard $J^{comp}$-invariant Riemannian metric (making it a Hermitian symmetric space). 
Then the pull-back of this metric to $\wt{M}$ is $\phi^{-1}(T)$-invariant. Recall that 
  $\overline{\Dt_{y_0}}\subset\Aut(\wt{M})$ fixes $y_0\in\wt{M}$, that $\phi(\wt{J})$
 is dense in $J$, and that $\phi(\overline{\Dt_{y_0}})$ is dense in $\overline{j'(\Dt_{y_0})}=T$.
Since the isotropy group at any point of a connected Riemannian manifold is compact (see for example 
\cite[II, Thm. 1.2]{Kob}),
we conclude that $\overline{\Dt_{y_0}}$ is compact. Hence $\phi(\overline{\Dt_{y_0}})=T$.
As $\Lie\,T$ spans all of $\Lie\,J$ as an $\Ad(j'(\wt{\Dt}))$-module over $\B{R}$,
the image of $\phi$ contains an open neighbourhood of the identity in $J$. Therefore $\phi:\wt{J}\to J$ is surjective. 
Using the fact that $\rho_0$ is a local isomorphism, we see 
that $\wt{J}$ is a Lie group and that $\phi$ is a topological covering.
\end{Emp}

\begin{Prop} \label{P:stab}
The stabilizer $J_y\subset J$ of each point $\rho_0(y)\in\rho_0(\wt{M})$ is compact.
\end{Prop}
\begin{pf}
We start from the following lemma.
\begin{Lem} \label{L:cart}
Every Cartan subgroup of $J$, contained in $J_y$, is compact.  
\end{Lem}
\begin{pf}
Suppose that the lemma is false, and let $C\subset J_y$ be a non-compact Cartan subgroup of $J$. 
Then the set 
$C':=\{c\in C:$ the group generated by $c$ is not relatively compact in $C\}$ is  a non-empty open 
subset of $C$. Set $\wt{C}:=\phi^{-1}(C)^0$. Then  $\wt{C}$ stabilizes $y\in\wt{M}$, and the subset 
$\wt{C}':=\{c\in\wt{C}:c$ is regular in 
$\wt{J}$ and $\phi(g)\in C'\}$ is a non-empty open subset of $\wt{C}$. Hence $\wt{J}':=\cup_{g\in\wt{J}} g\wt{C}'g^{-1}$ is an open non-empty subset of 
$\wt{J}$ (see for example the proof of \cite[Thm 7.101]{Kn}). 
It follows that there exists a $\dt$ belonging to $\wt{\Dt}\cap\wt{J}'$. Then by the very 
definition
of $\wt{J}'$, $\dt$ stabilizes some point on $\wt{M}$, and the subgroup, generated by $j'(\dt)$ is not
relatively compact in $J$. This contradicts to \rp{ell}.
\end{pf}
\begin{Cor} \label{C:form}
$\bold{J}$ is a real form of $(\PGL{d})^r$.
\end{Cor}
\begin{pf}
By our choice of $y_0$ and by \re{clos}, the group  $J_{y_0}$ always contains a Cartan subgroup of $J$. 
Hence the lemma shows that $J$ can not contain a factor of $\PGL{d}(\B{C})$. This shows the statement.
\end{pf}

Now we claim that every $J_y$ contains a Cartan subgroup of $J$. Notice that 
$J_y=J\cap\bold{P_y}(\B{C})\subset\bold{J}(\B{C})$ for some 
parabolic subgroup $\bold{P_y}\subset\bold{J_{\B{C}}}$. Hence  
$J_y=J\cap\bold{P_y}(\B{C})\cap\bold{\overline{P_y}}(\B{C})$, where $\bold{\overline{P_y}}$ is a parabolic of 
$\bold{J}_{\B{C}}$, complex conjugate to $\bold{P_y}$. The subgroup 
$\bold{P_y}\cap\bold{\overline{P_y}}\subset\bold{J}_{\B{C}}$ is rational over 
$\B{R}$, and, being an intersection of two parabolics, contains a maximal torus of $\bold{J}$ 
(see \cite[Ch. II, Cor of Thm. 5]{PR}). Therefore it contains a maximal $\B{R}$-rational torus, proving the claim.

Using the lemma, it now remains to show that $J_y$ is reductive. 
We have seen that $(\Lie\,J_y)_{\B{C}}$ contains a certain Cartan subalgebra $\C{T_y}:=(\Lie T_y)_{\B{C}}$ of 
$(\Lie\,J)_{\B{C}}$. Therefore it decomposes as a 
direct sum $\C{T_y}\oplus\oplus_{\beta\in\Phi}\C{J}^{\beta}$, where $\Phi$ is a subset of root system of 
 $(\Lie\,J)_{\B{C}}$ with respect to $\C{T_y}$, and $\C{J}^{\beta}$'s are the corresponding root spaces.
The Lie algebra $\Lie\,J$ is stable under the maximal compact subalgebra $\Lie\,T_y$ of $\C{T_y}$, and 
therefore
the complex conjugation  $\overline{ }$ of $(\Lie\,J)_{\B{C}}$ defined by its real form $\Lie\,J$ maps  
$\C{J}^{\beta}$ to $\C{J}^{-\beta}$. As  $(\Lie\,J_y)_{\B{C}}$ is  stable under this conjugation, we have
$\beta\in\Phi$ if and only if $-\beta\in\Phi$, so that $J_y$ is reductive.
\end{pf}


Take any $y\in\wt{M}$. Since $J_y$ is compact, we may assume that $J_y\subset J^{comp}$. 
Then $J_y$ is contained in the stabilizer $J^{comp}_y$ of $\rho_0(y)$ in $J^{comp}$. We know that 
$\dim\,\rho_0(\wt{M})=\dim\,((\B{P}^{d-1}_{\B{C}})^r)\an$, that $\dim\,J^0=\dim\,J^{comp}$, and that
$J^{comp}$ acts transitively on $((\B{P}^{d-1}_{\B{C}})^r)\an$. Therefore $J_y=J^{comp}_y=\bold{U_{d-1}}(\B{R})$, and 
the orbit $J(\rho_0(y))$ contains an open neighbourhood of $\rho_0(y)$. Since a point $y$ was chosen arbitrary and since 
$\rho_0(\wt{M})$ is connected, we see that $J_0$ acts transitively on $\rho_0(\wt{M})$.
Hence $\wt{J}$ acts transitively on $\wt{M}$, and $\rho_0:\wt{M}\to\rho_0(\wt{M})$ is a topological 
covering. As $J_y$ and $\rho_0(\wt{M})$ are connected, $J$ is connected as well.

Let $J_{y,i}\subset J_i$ be the $i$-th factor of $J_y$. 
As $\bo{U_{d-1}}(\B{R})$ is a maximal connected  proper
subgroup of $\PGU{d}(\B{R})$, we see by \cite[Ch. VIII, Thm. 6.1]{He} 
(considering separately compact and non-compact cases) that
$\rho_0(\wt{M})\cong J^0/J_y=\prod_{i=1}^r J_i/J_{y,i}$ is a Hermitian 
symmetric space, not containing Euclidean factors. In particular, (see \cite[Ch. VIII, Thm. 4.6]{He}) 
it is simply connected. Therefore
$\rho_0:\wt{M}\to\rho_0(\wt{M})$ and $\phi:\wt{J}\to J$ are isomorphisms. 

The kernel of $\wt{\pi}:\wt{\Dt}\to\Dt$ is isomorphic to the fundamental group of $M$. Hence it is 
normal in $\wt{\Dt}$ and discrete in $\Aut(\wt{M})$. Therefore it is central in $\wt{J}\cong J$, so that 
it is trivial. This shows that $\wt{\Dt}\cong\Dt$ and $\wt{M}\cong M$.

Recall that $M_S:=(\Dt\cap S)\bs M$ is a connected component of $(S\bs X_{\B{C}})\an$ for each 
$S\in\ff{\g}$, hence for a sufficiently small $S$ the projection $M\to M_S$
is a topological covering, $M_S$ is smooth, and its canonical bundle is ample 
(see \re{inv} and \rr{unif}). Therefore (see
\cite[IV, Thm. 5.3]{La}) $M_S$ and $M$ are measure-hyperbolic. In particular, 
$M$ has no compact factors. Then \cite[Ch. VIII, Thm. 7.13]{He} implies that
$((\B{P}^{d-1}_{\B{C}})^r)\an$ is a compact dual of $M$ and that $\rho_0$ is the Borel embedding. 
It follows that $M\cong\bbr$ and that $J\cong\Aut(\bbr)^0\cong\PGU{d-1,1}(\B{R})^r_{+}$.

{\bf Step 3}. Now we are going to determine $\Dt$. Recall that $\Dt$ is naturally embedded into
$\PGU{d-1,1}(\B{R})^r_+\times\g$ and that its projection to $\PGU{d-1,1}(\B{R})^r_+\subset\PGL{d}(\B{C})^r$
coincides with $j'$.

The following immediate generalization of \cite[Cor. 2.2.9]{V} plays the central role in this step. 

\begin{Lem} \label{L:rat}
 For each $\dt\in\Dt$ with an elliptic projection $\dt_{\be}\in\pgr^r_+$ there exists a representative 
$$\wt{\dt}=(\wt{\dt}_{\be},\wt{\dt}_v,\wt{f}_v,\wt{\dt}^{f;v})\in
\gr^r_{+}\times\prod_{i=1}^r\wt{D}_{w_i}\m\times\prod_{i=1}^r F_{v_i}\m\times \bold{G}(\afv)$$
such that the following hold:\par
a) Let us view $K$ as a subset of $K\otimes_{\B{Q}}\B{R}\cong
\B{C}^g$, of $\prod_{i=1}^r K_{w_i}$ and of $K\otimes_F\afv$ respectively. Then the characteristic 
polynomials of $\wt{\dt}_{\be},\wt{\dt}_v$ and 
$\wt{\dt}^{f;v}$ have their coefficients in $K$ and coincide.\par

b) $\wt{f}_v$ and the similitudes factor of $\wt{\dt}^{f,v}$ belong to $F\m$, viewed as a subset of  
$\prod_{i=1}^r F_{v_i}\m$ and of $(\afv)\m$ respectively, and coincide.
\end{Lem}

\begin{pf}
Since the arguments of the proof of \cite[Prop. 2.2.8 and Cor. 2.2.9]{V} work without changes in our
case, we omit the proof. Alternatively, the statement can be proved by very similar considerations to 
those of \rp{ell}, working with $P$ instead of $P'$.
\end{pf}

As in \cite[Prop. 1.6.1]{V}, $\Dt$ is a cocompact lattice of $\bo{J}(\B{R})^0\times\g$. 
Recall that the projection of $\Dt$ to $\bo{J}(\B{R})^0$ is injective and that the center $Z(\g)\cong
(\B{A}_K^f)\m/\overline{K}\m$ is compact. Hence the projection  $\Dt'$ of $\Dt$ to 
$\bo{J}(\B{R})^0\times P\g$ is also a cocompact lattice, isomorphic to $\Dt$. 
Now we proceed as in the proof of \cite[Thm. 2.2.5]{V}. 
For each finite set $N$ of finite places of $F$ and for each 
$S\in\ff{\B{A}_F^{f;N}}$ we define a cocompact lattice 
$\Dt^S\subset\prod_{i=1}^r\bo{J_i}(\B{R})\times\prod_{u\in N}\bo{PG\ii}(F_u)$
consisting of the projections of elements of $\Dt'$, whose components outside $N$ belongs to $S$ 
(see \cite[2.4.1]{V}). 

\begin{Lem} \label{L:irred}
The lattice $\Dt^S$ is irreducible (in the sense of \cite[Ch. III, Def. 5.9]{Ma}). 
\end{Lem}

\begin{pf}
We will prove a stronger assertion saying that for each $u\in\{\be_1,...,\be_r\}\cup N$ the projection 
$\pr_u$ of $\Dt'$ to the $u$'th component is injective. If not, set $\Dt^u:=\Ker\,\pr_u$.  
Then $\Dt^u$ is a non-trivial normal subgroup of $\Dt'$. Since the projection $\Dt'\to\bo{J}(\B{R})^0$ 
is dense and injective, the closure $\overline{\Dt^u_{\be}}$ of $\Dt^u_{\be}\subset\bo{J}(\B{R})^0$ 
is a normal non-trivial closed subgroup of $\bo{J}(\B{R})^0$. Therefore it is equal to 
$\prod_{i\in I}\bo{J_i}(\B{R})^0$ for some non-empty subset $I\subset\{1,...,r\}$.
Recall that the set of elliptic elements of $\pgr_+$ is open and non-empty (see \cite[Prop. 1.7.3]{V}).
Therefore there exists an element $\dt\in\Dt^u$, whose projection $\dt_{\be_i}\in\bo{J_i}(\B{R})^0$ is 
elliptic for all  $i\in I$ and trivial for all $i\notin I$. In particular, $\dt_{\be}$ 
fixes some point on $\bbr$. \rp{ell} and \rl{rat} now show that $\dt_u$ can not be equal to $1$, 
contradicting to our 
assumption. (Observe that the above proof is simpler than that of \cite[Prop. 2.4.5]{V}.)
\end{pf}

\begin{Emp} \label{E:Ma}
From this point all the arguments from \cite[2.3-2.6]{V} work without changes. First we show, as in 
\cite[2.4.5-2.4.6]{V}, that by Margulis' theorem \cite[Thm. (B), p. 298]{Ma} the subgroup $\Dt^S$ is 
arithmetic.
As in \cite[2.3 and 2.4.6]{V} (using \rp{ell}), we see that  an absolutely simple adjoint group 
$\bo{\wt{G}}$, defining $\Dt^S$, is defined over $F$ and does not depend on $N$ and $S$. We also get that 
$\Dt'\cong\Dt$ is naturally embedded into  $\bo{\wt{G}}(F)$. We conclude from \rl{rat}, as in 
\cite[2.5.2]{V}, that  $\bo{\wt{G}}$ is locally isomorphic to $\bo{PG\ii}$ at all places of $F$.
\end{Emp}

\begin{Emp} \label{E:inner}
Next we show directly that $\bo{\wt{G}}$ is an inner form of $\bo{PG\ii}$. In fact,
let $c$ be the element of $\text{H}^1(\Gal(\bar{F}/F),\Aut_{\bar{F}}(\bo{PG\ii}))$ corresponding to $\bo{\wt{G}}$, and
let $\C{D}$ be the Dynkin diagram of $\bo{PG\ii}/F$. Thus our aim is to show that the image $\bar{c}$ of 
$c$ in $\text{H}^1(\Gal(\bar{F}/F),\Aut(\C{D}))=\Hom(\Gal(\bar{F}/F),\B{Z}/2\B{Z})$ is trivial. 
Set $L$ be the subfield of $\bar{F}$ of fixed elements by $\Ker(\bar{c})$. 
Then $L$ is a quadratic or trivial extension of $F$, and
$\bar{c}\in\text{H}^1(\Gal(L/F),\Aut(\C{D}))$. Let $u$ be a finite prime of $F$, inert in $L$.
Then $\Gal(L/F)\cong\Gal(L_u/F_u)$. Hence it remains to show the triviality of the image 
$\bar{c}_u\in\text{H}^1(\Gal(L_u/F_u),\Aut(\C{D}))$ of $\bar{c}$. Since  $\bo{PG\ii}$ and $\bo{\wt{G}}$ are
locally isomorphic at all places of $F$, the image
$c_u\in\text{H}^1(\Gal(\bar{F}_u/F_u),\Aut_{\bar{F_u}}(\bo{PG\ii}))$ of $c$ is trivial.
But $\bar{c}_u$ is the image of $c_u$, and we are done.

It follows that  $\bo{\wt{G}}\cong\PGU{}(D',\al')$, where $(D',\al')$ is a form of $(D\ii,\al\ii)$, that is 
$D'$ is a central simple algebra over $K$ of dimension $d^2$, and $\al'$ is an involution of $D'$ of the second kind
over $F$. Thus we have avoided the use of the classification in \cite[2.5.1]{V}.
\end{Emp}

\begin{Emp}
As in \cite[2.5.4-2.5.6]{V}, we see that $\bo{\wt{G}}\cong\bo{PG\ii}$ and that the 
image of the natural embedding $\Dt\hra\g\cong\bo{G\ii}(\af)/\overline{\bo{Z(G\ii)}(F)}$ 
lies in  $\bo{\PG\ii}(F)_{+}$. Since $\Dt$ is cocompact in $\g$ and since $\bo{\PG\ii}(F)_{+}\subset\g$
is discrete, the image of $\Dt$ has a finite index in $\bo{\PG\ii}(F)_{+}$.

To show that these groups are isomorphic we can compare the volumes, as in \cite[2.6]{V}, using Kottwitz's 
results on Tamagawa measures \cite[Thm. 1]{Ko} and the Hirzebruch proportionality principle (\rp{Hirz}).
Alternatively one can argue as follows.

Let $h\in\bo{\PG\ii}(F)_{+}$ be an elliptic element, let $y\in\bbr$ be an elliptic point of $h$,
and let $\bo{T}$ be the maximal ($F$-rational) torus of $\bo{PG\ii}$ such that $\bo{T}(F)$ fixes $y$
and $h\in \bo{T}(F)$. Then by \rp{ell} there exist $x\in\prod_{i=1}^r\Omega^d_{K_{w_i}}(\B{C}_p)$ 
and $g\in P\g'$ such that $g(\Delta\cap\bo{T}(F))g^{-1}\subset \Gm_x\subset\bo{PG}(F)\subset P\g'$.
Since  $\bo{T}$ is connected, the subgroup $\bo{T}(F)\cap\Delta$ is Zariski dense in $\bo{T}$ 
(see \cite[Ch. V, Cor. 13.3]{Bo}).
Therefore the map $t\mapsto gtg^{-1}$ defines an algebraic
$F$-rational map from $\bo{T}$ to $\bo{PG}$. In particular,
$g h g^{-1}$ belongs to $\bo{PG}(F)$. Moreover, 
since the stabilizer of a point is an algebraic subgroup,
$g h g^{-1}\in\Gm_x$. \rp{ell} now shows that $h\in\Dt$. Since the set of elliptic
elements of $\bo{\PG\ii}(F)_{+}$ generates the whole group (see \cite[Ch. IX, Lem. 3.3]{Ma} and 
\cite[Prop. 1.7.3]{V}),
$\Dt\cong\bo{\PG\ii}(F)_{+}$.
\end{Emp}

 {\bf Step 4}. The proven above implies the existence of a 
$\g$-equivariant isomorphism $\wt{\varphi}:(X_{\B{C}})\an\isom(X\ii_{\B{C}})\an$. 
By GAGA results (compare \cite[Lem. 2.2.6]{V}) $\wt{\varphi}$ gives us an algebraic isomorphism 
$\varphi:X_{\B{C}}\isom X\ii_{\B{C}}$.

To see that  $\varphi$ is $E_v$-rational we will compare the action of the Galois group on the set of special points on the two sides.
Fix some special point $x\in X\ii(\B{C})$. Then $x$ corresponds to some maximal torus $\bo{T}\subset\bo{H\ii}$
and to some element $h'$ in the conjugacy class $M\ii$ of $h$. As in \cite[Lem. 3.1.5]{V}, $\bo{T}$
is equal to the intersection of $\bo{H\ii}$
with $\bo{R}_{L/\B{Q}}\B{G}_m$ for some maximal commutative subfield $L$ of $D\ii$, stable under 
$\al\ii$. Moreover, $L$ is a $CM$-field, and the restriction of $\al\ii$ to $L$ is the complex conjugation. For each $i=1,...,r$ the point $h'\in M\ii$ determines, as in \cite[3.1.6]{V}, an embedding
$\wt{\be}_i:L\hra\B{C}$, extending $\wt{\be}_i:K\hra\B{C}$.

The reflex field $\wt{E}\subset\B{C}$ of $(\bo{T},h)$ is the subfield  of the composite
$L_0$ of the fields $\wt{\be}_1(L),...,\wt{\be}_r(L)$ fixed elementwise by those automorphisms of $L_0/\B{Q}$
which permute the $\wt{\be}_i$'s.
Moreover, the canonical morphism $r':\bo{R}_{\wt{E}/\B{Q}}\B{G}_m\to\bo{T}$ 
(in the notation of \cite[3.1.4]{V}) is characterized as the unique morphism such that the composition map
$r'\circ N_{L_0/\wt{E}}:\bo{R}_{L_0/\B{Q}}\B{G}_m\to\bo{T}$ is given by 
$$r'(N_{L_0/\wt{E}}(l))=\prod_{i=1}^r\wt{\be}_i^{-1}(N_{L_0/\wt{\be}_i(L)}(\bar{l}/l))\in
\bo{T}(\B{Q})\subset L\m$$ 
for each $l\in L_0\m$, where by $\bar{l}$ we denote the complex conjugate of $l$ 
(compare \cite[3.1.7]{V}).

For each $i=1,...,r$ set $L_{w_i}:=L\otimes_K K_{w_i}$. Then $L_{w_i}$ is a maximal commutative subfield of
the skew field $D\ii\otimes_K K_{w_i}\cong\wt{D}_{w_i}$. As in \re{p-ad}, each $\wt{\be}_i$ defines a continuous
embedding $\wt{\al}_i:L_{w_i}\hra\B{C}_p$. Moreover, the composite field of the $\wt{\al}_i(L)$'s is the completion 
$\wt{E}_v$ of $\wt{E}\subset\B{C}_p$ (compare \rl{comp}).

For each $l\in\wt{E}\m_v$ define an element 
$$\wt{r}(l):=(\wt{\al}_1^{-1}(N_{\wt{E}_v/\wt{\al}_1(L_{w_1})}(l^{-1})),
...,\wt{\al}_r^{-1}(N_{\wt{E}_v/\wt{\al}_r(L_{w_r})}(l^{-1}));1)$$
in $\prod_{i=1}^r L_{w_i}\m\times\g'\subset\prod_{i=1}^r \wt{D}_{w_i}\m\times\g'=\wt{\g}$. Let $\theta_{\wt{E}_v}:\wt{E}\m_v\to\Gal(\wt{E}_v^{ab}/\wt{E}_v)$
be the Artin homomorphism, sending the uniformizer to the arithmetic Frobenius automorphism.
By the definition of $X\ii$
(see for example \cite[3.1.4, 3.1.12]{V}), $x$ is rational over $\wt{E}_v^{ab}$ and for each
$l\in\wt{E}\m_v$ the element $\theta_{\wt{E}_v}(l)\in\Gal(\wt{E}_v^{ab}/\wt{E}_v)$ maps $x$ to $\wt{r}(l)(x)$
(compare \cite[Prop. 3.1.9]{V}). 
 Since these properties (for all special points) characterize the weakly-canonical models,
it will suffice to show that

i) the point $y:=\varphi^{-1}(x)\in X(\B{C}_p)$ is rational over $\wt{E}_v^{ab}$;

ii) $\theta_{\wt{E}_v}(l)(y)=\wt{r}(l)(y)$ for each $l\in\wt{E}\m_v$. 

Let $[\wt{y},g]\in\SI\times\g'$ be a representative of $y$. Then  $[\sigma(\wt{y}),g]$  is 
a representative of $\sigma(y)$ for each  (not necessary continuous) $\sigma\in\Aut(\B{C}_p/\wt{E}_v)$.
As in \cite[Prop. 3.2.3]{V}, we show that for each $i=1,...,r$ there exists an embedding 
$L_{w_i}\hra\Mat_d(K_{w_i})$ such that $\wt{y}$ belongs to the image of the corresponding
embedding $\prod_{i=1}^r\Sigma_{L_{w_i}}^1\otimes_{K_{w_i}}E_v\hra\SI$ (see \cite[1.4.4]{V}). 
Then the statement follows from \cite[Lem. 1.4.3]{V}.
This completes the proof of the First Main Theorem.

\begin{Rem} \label{R:comp}
a) Let $(X\ii)'$ be the weakly-canonical model over $E_v$ of the Shimura variety, corresponding to the same group
$\bold{H\ii}$ as $X\ii$ and to a homomorphism  $h':\bo{R}_{\B{C}/\B{R}}\B{G}_m\to\bo{H\ii}_{\B{R}}$ defined by 
requirement that for  each $z\in\B{C}\m$ we have
$$h(z)=(\underbrace{\diag(\bar{z},...,\bar{z},z)^{-1};...;\diag(\bar{z},...,\bar{z},z)^{-1}}_{r};\underbrace{I_d;...;I_d}_{g-r})\in\GU{d-1,1}(\B{R})^r\times\GU{d}(\B{R})^{g-r}.$$
Then $(X\ii)'$ is a certain abelian unramified twist of $X\ii$. Moreover, comparing the action of the Galois group 
on the set of special points, one can check that the difference of the Galois actions on $X\ii$ and on 
$(X\ii)'$ is given by some explicit homomorphism $\Aut(\B{C}/E_v)\to\g_1:=\{1\}\times\prod_{i=1}^r F_{v_i}\m\times\{1\}
\subset\wt{\g}=\bo{G\ii}(\af)$. In particular, the quotients $\g_1\bs X\ii$ and $\g_1\bs(X\ii)'$ (see \rl{quot2})
are naturally isomorphic (over $E_v$).

b) Rapoport and Zink work with $(X\ii)'$ instead of $X\ii$, therefore their Shimura varieties are uniformized 
by the corresponding twist (see \cite[Prop. 6.49]{RZ2}) of the product of Drinfeld upper half-spaces, and 
\cite[Thm. 6.50]{RZ2} follows from our First Main Theorem by twisting.
\end{Rem}

 {\bf Step 5}. We have already shown in the second and the third steps that the homomorphism 
$j':\bo{PG\ii}(F)_+\to\PGL{d}(\B{C})^r\cong\prod_{i=1}^r\bo{PG\ii}(K_{\wt{\be}_i})$ is a diagonal embedding.
Arguing as in \cite[4.6.2-4.6.4 and 4.5.5]{V} we show that this holds for the whole $j$.
For the descent to $E_v$ we use the uniqueness arguments of \cite[4.7]{V}. This completes the proof of the 
Second Main Theorem.

\begin{Cor}  \label{C:cover}
 Let $X_{\Gm}$ be a projective scheme over a $p$-adic field $L$ such that $X\an_{\Gm}\cong\Gm\bs\prod_{i=1}^r
\Omega^d_{L_i}\otimes_{L_i} L$ for $p$-adic subfields $L_1,...,L_r$ of $L$, a natural number $d\geq 2$,
and a  torsion-free cocompact lattice $\Gm\subset\prod_{i=1}^r\PGL{d}(L_i)$ (see \rl{quot}).
Then for each field embedding $L\hra\B{C}$
the universal cover of $(X_{\Gm,\B{C}})\an$ is $\bbr$.
\end{Cor}

\begin{pf}
For curves ($r=1$ and $d=2$) the statement follows from the fact that $X_{\Gm}$ is a smooth curve of 
genus greater then $1$ (see \cite[Thm. 3.3]{Mum}).

In the remaining cases, the rank of $\prod_{i=1}^r\PGL{d}(K_{w_i})$ is $r(d-1)\geq 2$. If $\Gm$ is irreducible, then 
it is arithmetic by \cite[Thm. (A), p. 298]{Ma}.
Therefore there exists a number field $F$, an adjoint absolutely simple group $\bo{G'}$
over $F$, and finite primes $v_1,...,v_r$ of $F$
such that $\bo{G'}(F_{v_i})\cong\PGL{d}(L_i)$ for each $i=1,...,r$, the group $\bo{G'}(F_{\be_i})$ is compact
for each archimedean completion $F_{\be_i}$ of $F$, and $\Gm$ is commensurable with 
$\bo{G'}(F)\cap S\subset\prod_{i=1}^r\bo{G'}(F_{v_i})\cong\prod_{i=1}^r\PGL{d}(K_{w_i})$ for every 
$S\in\ff{\bo{G'}(\afv)}$.
Since the universal cover of $(X_{\Gm,\B{C}})\an$ does not change after  replacement of $\Gm$ by another 
torsion-free subgroup of $\prod_{i=1}^r\PGL{d}(K_{w_i})$, commensurable with $\Gm$, we may assume that
$\Gm=\bo{G'}(F)\cap S$ for a sufficiently small $S$.
Now we can apply the first two steps of the proof of the  Main Theorems (using $\bo{G'}$ instead of 
$\bo{G}$ to conclude that $(X_{\Gm,\B{C}})\an\cong \Dt_S\bs\bbr$ for some torsion-free
cocompact lattice $\Dt_S\subset\pgr_+^r$. This proves the statement for irreducible $\Gm$'s.

If $\Gm$ is reducible, then possibly replacing $\Gm$ by a subgroup of finite index
we may assume that there exists a non-trivial  proper subset $I\subset\{1,...,r\}$ and torsion-free
cocompact lattices $\Gm_1\subset\prod_{i\in I}\PGL{d}(K_{w_i})$ and $\Gm_2\subset\prod_{i\notin I}\PGL{d}(K_{w_i})$
such that $\Gm=\Gm_1\times\Gm_2$. Then $\Gm\bs\prod_{i=1}^r
\Omega^d_{L_i}\otimes_{L_i} L\cong(\Gm_1\bs\prod_{i\in I}\Omega^d_{L_i}\otimes_{L_i} L)\times
(\Gm_2\bs\prod_{i\notin I}\Omega^d_{L_i}\otimes_{L_i} L)$, and the statement follows from the previous cases
by induction on the number of factors.
\end{pf}

\begin{Rem}
Observe that the proof uses neither the classification of algebraic groups nor of
Hermitian symmetric domains.
\end{Rem}

\begin{Rem} \label{R:Yau}
We know that the canonical bundle of $X_{\Gm}$ is ample (\rl{quot}) and that the Chern numbers of  $X_{\Gm}$ are proportional to those of $(\B{P}^{d-1})^r$ (see \rp{Hirz}). Therefore in the case $r=1$ the 
corollary follows immediately from Yau's theorem (\cite{Ya1}).
\end{Rem}

\begin{Rem} \label{R:equiv}
The proof shows that for arithmetic $\Gm's$ the corollary is a direct consequence of the results of the Step 2.
Conversely, the arguments of \cite[proof of Prop. 1.6.1]{V} show that the corollary implies that
$M\cong (B^{d-1})^r$. In particular, in the case $r=1$ the differential geometric part of the 
First Main Theorem immediately follows from Yau's Theorem. 
\end{Rem}

\section{Another approach}

In this section we give another approach to the differential-geometric part of the proof
of the First Main Theorem. Using Yau's theorem on the existence and the uniqueness of the 
K\"ahler-Einstein metric, (but not appealing to principal bundles and to the structure theory of reductive groups)
we will show directly that (in the notation of Step 1 of Section 3) 
$M$ is a Hermitian symmetric domain. Here we use the ideas from Milne's unpublished manuscript \cite{Mi2}.

%
%
%


\begin{Lem} \label{L:metric}
There exists a $\Dt$-invariant K\"ahler metric $g$ on $M$.
\end{Lem}

\begin{pf}
Recall the following theorem of Yau (see \cite[Thm. 5]{Ya}). Let $V$ be a smooth complex projective 
variety whose canonical bundle is ample.
Then there exists a unique K\"{a}hler metric $g_V$ on $V\an$, called the K\"ahler-Einstein metric, 
whose Ricci tensor satisfies Ric$(g_V)=-g_V$ (notice that  the sign of Yau's Ricci tensor differs from the 
standard one).

We are going to apply this theorem to $M_S:=(\Dt\cap S)\bs M$ for sufficiently small $S\in\ff{\wt{\g}}$
(see the last paragraph of Step 2 of Section 3). 

Let $g_S$ be the K\"ahler-Einstein metric on $M_S$, and let $g_{(S)}$ be its inverse image on $M$.
Since K\"ahler-Einstein metric is unique we conclude that
$g:=g_{(S)}$ does not depend on $S$ and that it is  $\Dt$-invariant.

\end{pf}


Since $(M,g)$ is a K\"ahler manifold, the group Is$(M,g)$ of its isometries is a Lie group
(see \cite[II, Thm. 1.2]{Kob}). Since $g$  is $\Dt$-invariant, 
$\Dt$ is contained in  Is$(M,g)$, therefore its closure $G$ is a Lie subgroup of  Is$(M,g)$.
Notice also that $G$ is also contained in $\Aut(M)$. For each $y\in M$ we denote the stabilizers of $y$ 
in $\Dt$ and in $G$ by $\Dt_y$ and $G_y$ respectively.


\begin{Lem} \label{L:mult}
There exist a point $y\in M$ and an element $\dt\in G_y$ of order $4$, inducing the multiplication 
by $\sqrt{-1}$ on $T_y(M)$.
\end{Lem}

\begin{pf}
(Compare the proof of \rp{ell}). Let $x\in\prod_{i=1}^r\Omega^d_{K_{w_i}}(\B{C}_p)\subset\B{P}^{d-1}(\B{C}_p)^r$
be an elliptic point for the action of $\PG(F)$. Then the point $[x,1]\in[\prod_{i=1}^r\Omega^d_{K_{w_i}}
\otimes_{K_{w_i}}E_v\times (P\g')^{disc}]/\PG(F)\cong (X')\an$ corresponds to a certain point
$z\in X'(\B{C}_p)=X'(\B{C})$ and thethefore to a point $[y,g]\in [M\times(P\g')^{disc}]/\Dt\cong(X'_{\B{C}})\an$. Therefore the map
$\dt\mapsto g\dt g^{-1}$ defines an isomorphism between $\Dt_y$ and the stabilizer 
$\Gm_x\subset\PG(F)\subset P\g'$ of $x$. Now we identify $\B{C}=\B{C}_p$-vector spaces
$T_y(M)$ and $T_x(\B{P}^{d-1}_{\B{C}_p})^r=T_x(\B{P}^{d-1}_{\B{C}})^r$ via the sequence of canonical 
isomorphisms 
$$T_y(M)\cong T_{[y,1]}(X'_{\B{C}})\an\cong T_{[y,g]}(X'_{\B{C}})\an\cong T_z(X'_{\B{C}}=T_z(X'_{\B{C}_p}\cong$$ 
$$\cong T_{[x,1]}(X'_{\B{C}_p})\an\cong T_x(\prod_{i=1}^r\Omega^d_{K_{w_i}}\otimes_{K_{w_i}}\B{C}_p)
\cong T_x(\B{P}^{d-1}_{\B{C}_p})^r.$$
Then the above isomorphism $\Dt_y\isom\Gm_x$ commutes with embeddings $\Dt_y\hra\Aut(T_y(M))$
and $\Gm_x\hra\Aut(T_x(\B{P}^{d-1}_{\B{C}})^r)$, hence it extends by continuity to the isomorphism between 
the closures $\overline{\Dt_y}\subset\Aut(T_y(M))$ and $\overline{\Gm_x}\subset\Aut(T_x(\B{P}^{d-1}_{\B{C}})^r)$. Since these closures are naturally isomorphic to the closures
in $G$ and in $\PH(\B{R})$ respectively, it remains to show the existence of an element in
$\overline{\Gm_x}\subset\PH(\B{R})\cong\PGU{d}(\B{R})^r$, which has order $4$ and induces the multiplication 
by $\sqrt{-1}$ on $T_x(\B{P}^{d-1}_{\B{C}})^r$. As $\overline{\Gm_x}$ is a maximal compact torus
of $\PGU{d}(\B{R})^r$ (see the last paragraph of the proof of \rp{ell}), this is clear.
\end{pf}

\begin{Prop} \label{P:Herm}
$M$ is a Hermitian symmetric domain.
\end{Prop}

\begin{pf}
First we show that $G$ acts transitively on $M$. Let $y\in M$ be as in the lemma,
and let $N$ be the orbit $Gy$. Consider the map $\varphi_y:$ Is$(M,g)\to M;\;\dt\mapsto \dt(y)$.
By \cite[p.41]{Kob}, it factors as a composition of a closed embedding of  Is$(M,g)$
into the bundle  $O(M)$ of orthogonal frames over $M$ and the projection of $O(M)$ to $M$.
This shows that $\varphi_y$ is proper, therefore
$N$ is closed in $M$. By \cite[II, Thm. 3.2 and 4.2]{He},  
$N\cong G/G_y$ then has a canonical structure of a $G$-equivariant closed real submanifold of $M$.
In particular, its tangent space $T_y(N)\subset T_y(M)$ is invariant under the action of $G_y$.
By \rl{mult} it then invariant under the multiplication by $\sqrt{-1}$.
Since the group $G$ acts on $N$ transitively, $N$ is an almost complex submanifold of the complex
manifold $M$. Therefore $N$ has a unique $G$-invariant structure of a closed complex submanifold
of $M$ (see \cite[Ch. VIII, Thm. 1.2 and p.284]{He}).
Therefore $\wt{N}:=[N\times\g^{disc}]/\Dt$ is a $\g$-invariant closed (complex) analytic subspace of 
$(X_{\B{C}})\an\cong[M\times\g^{disc}]/\Dt$.

\begin{Lem} \label{L:alg}
There exists a closed $\g$-invariant subscheme $Y$ of $X_{\B{C}}$ 
such that $Y\an\cong\wt{N}$.
\end{Lem}

\begin{pf}
Take an $S\in\ff{\g}$. Then the quotient map $\pi_S:(X_{\B{C}})\an\to (S\bs X_{\B{C}})\an$ is open
and $\wt{N}$ is  $S$-invariant.
Therefore the quotient $\pi_S(\wt{N})\cong S\bs[N\times\g^{disc}]/\Dt$ is a closed 
analytic subspace  of $(S\bs X_{\B{C}})\an$. The scheme $S\bs X_{\B{C}}$
is projective, hence by GAGA (see for example \cite[Cor. 1.2.2]{V}) there exists a projective subscheme
$Y_S$ of $S\bs X_{\B{C}}$ such that $Y_S\an\cong S\bs[N\times\g^{disc}]/\Dt$.
The inverse limit of the $Y_S$'s satisfies the conditions of the lemma.
\end{pf}

Since every $\g$-orbit in $X_{\B{C}}$ is Zariski dense (as in \cite[Prop. 1.5.3 e)]{V}), 
we obtain that $Y=X_{\B{C}}$. Hence $N=M$, so that $G$ acts transitively on $M$, as claimed.

Let $\dt\in G$ be as in \rl{mult}. Then for each $z=g(y)\in M$  with $g\in G$ we have an involutive 
holomorphic isometry $g\dt^2 g^{-1}$ of $M$ with an isolated fixed point $z$. This means that
$M$ is a Hermitian symmetric space. Moreover, as in Step 2 of Section 3 we see that 
$M$ is actually a Hermitian symmetric domain. 
\end{pf}

\begin{Rem} \label{R:alter}
Knowing the proposition one can obtain an alternative proof of the results of Step 2 of Section 3.
For shoing this one can simply procceed along the lines of the above mentioned step, 
which would be much easier assuming the proposition. Alternatively,
one can apply Margulis' superrigidity theorem (as in \cite[4.6]{V}). 
\end{Rem}

\section{The case of quaternion Shimura varieties} \label{S:quat}

In this section we deduce a $p$-adic uniformization of quaternion Shimura varieties from the unitary case
(compare \cite{BZ}, where the moduli approach (see \cite[\S6]{De}) is used).

\begin{Emp}
Let $F$ be a totally real number field of degree $g$ over $\B{Q}$,  let 
$\be_1,...,\be_g$ be the archimedean completions of $F$, and  let $1\leq r\leq g$ be a natural number.
 Let $D\ii$ be a quaternion algebra over $F$, split at $\be_1,...,\be_r$ and ramified at 
$\be_{r+1},...,\be_g$. Let $\bo{G\ii}=\GL{1}(D\ii)$ be the algebraic group over $F$, 
characterized by  $\bo{G\ii}(R)=(D\ii\otimes_F R)\m$ for each $F$-algebra $R$, and set  
$\bo{H\ii}:=\bo{R}_{F/\B{Q}}\bo{G\ii}$. 

Define a homomorphism $h:\bo{R}_{\B{C}/\B{R}}\B{G}_m\to\bo{H\ii}_{\B{R}}$ by
requiring for each $z=x+iy\in\B{C}\m$
$$h(z)=(\underbrace{{\begin{pmatrix}
                    x & y \\ -y & x
                    \end{pmatrix}}^{-1};...;
                   {\begin{pmatrix}
                    x & y \\ -y & x
                    \end{pmatrix}}^{-1}}_{r};\underbrace{\begin{pmatrix}
                    1 & 0 \\ 0 & 1
                    \end{pmatrix};...;
                   \begin{pmatrix}
                    1 & 0 \\ 0 & 1
                    \end{pmatrix}}_{g-r})\in \GL{2}(\B{R})^r\times(\B{H}\m)^{g-r}\cong\bo{H\ii}(\B{R}),$$
where $\B{H}$ denotes the field of Hamilton's quaternions.
Then  $M\ii$, the $\bo{H\ii}(\B{R})$-conjugacy class of $h$,  is isomorphic to $\h^r$.
The reflex field $E\subset\B{C}$ of $(\bo{H\ii},M\ii)$ is the subfield of the composite $F_0$
of the fields $\be_1(F),...,\be_r(F)$ fixed elementwise by those automorphisms of $F_0/\B{Q}$,
which permute the $\be_i$'s ($i=1,...,r$).

Let $v$ be a finite prime of $E$, and let $v_1,...,v_r$ be the finite primes of $F$, corresponding to $\be_1,...,\be_r$
as in \re{p-ad}. Assume that the $v_i$'s are distinct and that $D\ii$ is ramified at $v_1,...,v_r$.
Then, as in \rl{comp}, we see that for each $i=1,...,r$ we have a natural embedding $\al_i:F_{v_i}\hra E_v$ and that
$E_v$ is the composite field of the $\al_i(F_{v_i})$'s. Assume that 
the $v_i$'s are distinct and that $D\ii$ is ramified at $v_1,...,v_r$.

Let ${X}\ii$ be the weakly-canonical model over $E_v$ of the Shimura variety corresponding to  $(\bo{H\ii}, M\ii)$.
Let $\check{M}\ii$ be the Grassmann variety corresponding to the pair $(\bo{H\ii}, M\ii)$, and let $W\ii$ be a $\bo{PH\ii_{E_v}}$-equivariant vector bundle on $\check{M}\ii_{E_v}$. Let $P\ii$ (resp. $V\ii=V\ii(W\ii)$) be the 
canonical model of the standard $\bo{PH\ii_{E_v}}$-principal
bundle over $X\ii$ (resp. the automorphic vector bundle on $X\ii$ corresponding to $W\ii$).
\end{Emp}

\begin{Emp}   
Let $D$ be the quaternion algebra over $F$ ramified at $\be_1,...,\be_r$, split at 
$v_1,...,v_r$, and locally isomorphic to $D\ii$ at all other places of $F$. Notice that $D$ is definite.

Define algebraic groups $\bo{G}$ and $\bo{H}$ by  $\bo{G}:=\GL{1}(D)$ and $\bo{H}:=\bo{R}_{F/\B{Q}}\bo{G}$. 
For each $i=1,...,r$ fix a quaternion  division 
algebra $\wt{D}_{v_i}$ over $F_{v_i}$. Set $\g':=\bo{G}(\afv)$ and $\wt{\g}:=\prod_{i=1}^r \wt{D}_{v_i}\m\times\g'$.
Fix algebra isomorphisms  $D\otimes_F F_{v_i}\cong\Mat_2(F_{v_i}),\quad D\ii\otimes_F F_{v_i}\cong 
\wt{D}_{v_i}$ and  $D\otimes_F\afv\cong D\ii\otimes_F\afv$. They induce isomorphisms 
$\bo{G}(\af)\cong\prod_{i=1}^r\GL{2}(F_{v_i})\times \g'$ and $\bo{G\ii}(\af)\cong\wt{\g}$.
We will often identify the groups  $\bo{G\ii}(\af)$ and $\wt{\g}$ by means of the last isomorphism.

Let $X$ be a scheme over $E_v$ with a $\wt{\g}$-action such that for each
$S\in\ff{\wt{\g}}$ the quotient $S\bs X$ is projective and satisfies 
$$(S\bs X)\an\cong S\bs [\prod_{i=1}^r\Sigma^2_{F_{v_i}}\otimes_{F_{v_i}}E_v\times\g']/\bo{G}(F)$$
(compare \rp{alg} and \cite[Prop. 1.5.2]{V}). Let $P$ (resp. $V$) be the $\bo{PH_{E_v}}$-principal 
(resp. the automorphic vector) bundle over $X$ defined as in \re{tw}.
\end{Emp}

\begin{Thm} \label{T:quat}
There exists a $\wt{\g}$-equivariant isomorphism $\varphi:X\isom X\ii$. Moreover,  $\varphi$ lifts
to $\wt{\g}$-equivariant isomorphisms

a) $\varphi_V:V\isom V\ii$ of automorphic vector bundles and

b) $\varphi_P:P\tw\isom P\ii$ of $\bo{PH\ii_{E_v}}$-principal bundles.
\end{Thm}
 
\begin{pf}
The theorem can be proved exactly by the same argument as the Main Theorems. In fact, the proof would be
much easier technically. 
Instead we will embed our schemes into certain quotients of the corresponding schemes in the unitary case and 
will deduce the theorem from the Main Theorems. 

First we define the corresponding objects in the 
unitary case. Let $K$ be a totally imaginary quadratic extension of $F$ in which all $v_1,...,v_r$ split. Then for each 
$i=1,...,r$ we can and do identify $F_{v_i}$ (resp. $\wt{D}_{v_i})$ with $K_{w_i}$ (resp. $\wt{D}_{w_i}$),
defined as in Section 2.
As a consequence, the field $E_v$ has the same meaning, as in Section 2.
Define the involution of the second kind $\al$ (resp. $\al\ii$) on $D\otimes_F K$ (resp.  $D\ii\otimes_F K$),
as the tensor product of the main involution of $D$ (resp. $D\ii$) and the conjugation of $K$ over $F$.

Without further mention we will write all objects in the unitary case, defined from the above data as in Section 2, 
with an upper subscript $ ^{un}$, 
meaning ''unitary''. For example, we write $\bo{G\un}$ for $\GU{}(D\otimes_F K,\al)$ and $\bo{(G\ii)\un}$ for 
$\GU{}(D\ii\otimes_F K,\al\ii)$. Then the above algebra isomorphisms identify, as in \re{not2}, 
$\wt{\g}\un$ with $\bo{(G\ii)\un}(\af)$. Define a subgroup $\g_1\subset\wt{\g}\un=\bo{(G\ii)\un}(\af)$ as in \rr{comp}.

The diagonal embeddings $\bo{G}\hra\bo{G\un}$ and $\bo{G\ii}\hra\bo{(G\ii)\un}$ define us
embeddings
$$\wt{\g}=\prod_{i=1}^r \wt{D}_{v_i}\m\times\bo{G}(\afv)\hra\prod_{i=1}^r \wt{D}_{w_i}\m\times\bo{G\un}(\afv)=
\g_1\bs\wt{\g}\un\text{ and }$$
$$\bo{G\ii}(\af)\hra\bo{(G\ii)\un}(\af)\surj\g_1\bs\bo{(G\ii)\un}(\af),$$
commuting with the above identifications.

 The following lemma and its obvious analogs will assure the existence of the  quotients, used below.
\begin{Lem} \label{L:quot2}
For each closed subgroup $H\subset Z(\wt{\g}\un)$ the quotient scheme $H\bs X\un$ exists.
\end{Lem}

\begin{pf}
(ompare \cite[Lem. 1.3.11 c) and Prop. 1.5.3 c)]{V}). In the notation of  Step 1 of Section 3, the group $Z(\g\un)=Z(\wt{\g}\un)/\g_0\un$ is compact. 
Hence for each $S\in\ff{\wt{\g}\un}$ the subgroup $H\cdot S$ of $\wt{\g}\un$ is open
and compact modulo $\g_0\un$. Since $\g_0\un$ acts trivially on $X\un$,
each quotient $(H\cdot S)\bs X\un$ exists, and their inverse limit is $H\bs X\un$.
\end{pf}

Next we notice that the natural equivariant embeddings 
$$[\prod_{i=1}^r\Sigma_{F_{v_i}}^2\times\bo{G}(\afv)]/\bo{G}(F)\hra
[\prod_{i=1}^r\Sigma_{K_{w_i}}^2\times\bo{G\un}(\afv)]/\bo{G\un}(F)\text{  and  } $$
$$[\h^r\times\bo{G\ii}(\af)]/\bo{G\ii}(F)\hra [\h^r\times(\g_1\bs\bo{(G\ii)\un}(\af))]/\bo{(G\ii)\un}(F)$$
define by GAGA equivariant embeddings $X\hra\g_1\bs X\un$ and $i:X\ii_{\B{C}}\hra\g_1\bs (X\ii)\un_{\B{C}}$.

To show that $i$ is  $E_v$-rational we observe that the natural embedding $\bo{H\ii}\hra\bo{(H\ii)\un}$ maps homomorphism
$h$ to the conjugate of $(h\un)'$ (in the notation of \rr{comp}). Therefore the natural embedding 
$X\ii\hra((X\ii)\un)'$ is $E_v$-rational (see \cite[Cor. 5.4]{De}). Dividing by $\g_1$ we get our statement from  
\rr{comp}. (Alternatively one can check it directly by calculating the action of the Galois group on the set of 
special points (as in Step 4 of Section 3, where the unitary case is treated).

Next we observe that after dividing by the center we have $\bo{PG}\cong\bo{PG\un}$ and $\bo{PG\ii}\cong\bo{P(G\ii)\un}$,
implying that the above embeddings of schemes induce isomorphisms 
$Z(\wt{\g})\bs X\isom Z(\wt{\g}\un)\bs X\un$ and $Z(\bo{G\ii}(\af))\bs X\ii\isom Z(\bo{(G\ii)\un}(\af))\bs (X\ii)\un$.

Let $\varphi\un:X\un\isom (X\ii)\un$ be a $\wt{\g}\un$-equivariant isomorphism as in the First Main Theorem.
It gives us a  $\g_1\bs\wt{\g}\un$-equivariant isomorphism 
$\bar{\varphi}:=\g_1\bs\varphi\un:\g_1\bs X\un\isom\g_1\bs (X\ii)\un$. Fix an $x\in X\subset\g_1\bs X\un$.
As it was explained above, there exists $g\in Z(\wt{\g}\un)$ such that 
$g\bar{\varphi}(x)\in X\ii\subset\g_1\bs (X\ii)\un$. Hence $g\bar{\varphi}$ maps the $\wt{\g}$-orbit of $x$ into $X\ii$.
Since each $\bo{G\ii}(\af)(=\wt{\g})$-orbit is Zariski dense in both $X$ and $X\ii$ (see
\cite[Prop. 1.3.8 and 1.5.3]{V}), $g\bar{\varphi}$ restricts to a  
$\wt{\g}$-equivariant isomorphism $\varphi:X\isom X\ii$.
The existence of the liftings $\varphi_V$ and $\varphi_P$ is now an immediate consequence of the Second Main Theorem. 
\end{pf}

\section{On the generality of our results} \label{S:gen}

In this section we will show that the Shimura varieties treated in this paper are the most general one (up to a central modification)
which can be $p$-adically uniformized by the product of Drinfeld upper half-spaces. 

Let $\bo{\wt{H}}$ be a reductive group over $\B{Q}$, and let $\wt{M}$ be a conjugacy class of homomorphisms
$\bo{R}_{\B{C}/\B{R}}\B{G}_m\to\bo{\wt{H}}_{\B{R}}$ satisfying Deligne's axioms. Let $E'\subset\B{C}$ 
be the reflex field of $(\bo{\wt{H}},\wt{M})$ and let $v'$ be a finite prime of $E'$. Let $X'$ be the 
canonical model of the Shimura variety corresponding to  $(\bo{\wt{H}},\wt{M})$. Suppose that there exists
$T\in\ff{\bo{\wt{H}}(\B{A}^f)}$ such that the projection $X'\to T\bs X'$ is \'etale and that some connected
component of $( T\bs X'\otimes_{E'}E'_{v'})\an$ is a finite \'etale cover of an analytic space of the form
$\Gm\bs\prod_{i=1}^r\Omega^d_{E_i}\otimes_{E_i}E'_{v'}$ for some $p$-adic subfields $E_1,...,E_r$ of $E'_{v'}$,
some natural number $d\geq 2$, 
and some irreducible arithmetic torsion-free cocompact lattice $\Gm\subset\prod_{i=1}^r\PGL{d}(E_i)$
(compare \rr{unif}). Recall that when $r(d-1)\geq 2$, that is in all cases except for curves,
the arithmeticity follows from the irreducibility.

\begin{Cl} \label{C:isom}
$\bo{P\wt{H}}$ is isomorphic to the group $\bo{R}_{F/\B{Q}}\PGU{}(D\ii,\al\ii)$ for some $F$, $D\ii$ and $\al\ii$
satisfying the assumptions of  \rs{stat}.
\end{Cl}

\begin{pf}
Since $\Gm$ is arithmetic, there exist $F$, $\bo{G'}$ and $v_1,...,v_r$ as in the proof of \rco{cover}.
In particular, $\bo{G'}$ is a form of $\PGL{d}$, $F$ is totally real, and each $\bo{G'_{F_{\be_i}}}$ is 
isomorphic to $\PGU{d}$, the unique $\B{R}$-unisotropic form of $\PGL{d}$. 

Moreover, by the classification of simple algebraic groups (see \cite{Ti}), $\bo{G'}\cong\PGU{}(D,\al)$ for some central simple algebra $D$ of dimension $d^2$
over a totally imaginary quadratic extension $K$ over $F$ and some involution of the second kind $\al$
of $D$ over $F$. Since  $\bo{G'_{F_{v_i}}}\cong\PGL{d}$ for each $i=1,...,r$, the $v_i$'s split in $K$, and $D$ 
splits at each prime $w_i$ of $K$ lying over $v_i$. Hence $D$ and $\al$ satisfy the assumptions 
of \rs{stat}. 

Let $D\ii$ and $\al\ii$ correspond to $D$ and $\al$, as in \rs{stat}, and set $\bo{(G')\ii}:=\PGU{}(D\ii,\al\ii)$.
Our assumptions together with the First Main Theorem imply that 
 $\bbr$ is the universal cover
of  some connected component $\wt{X}_T^0$ of 
$(T\bs X'\otimes_{E'}\B{C})\an$. Furthermore,  the fundamental group $\Gm_1$ of  $\wt{X}_T^0$ is commensurable with 
$\Gm_2:=\bo{(G')\ii}(F)_+\cap S\subset\PGU{d-1,1}(\B{R})^r$ for every $S\in\ff{\bo{(G')\ii}(\af)}$.
In particular, the arithmetic lattice $\Gm_1\subset \bo{P\wt{H}}(\B{R})$ is irreducible.
By Deligne's assumption, $\bo{P\wt{H}}$ has no $\B{Q}$-rational $\B{R}$-unisotropic factors, so it is 
$\B{Q}$-simple. Hence there exists a number field $F'$ and an  adjoint absolutely simple group
$\bo{\wt{G}}$ over $F'$ such that $\bo{P\wt{H}}\cong\bo{R}_{F'/\B{Q}}\bo{\wt{G}}$. 

Notice now that $\Gm_1$ and $\Gm_2$ are Zariski dense in $\bo{\wt{G}}$ and $\bo{(G')\ii}$ respectively 
(see \cite[I, Prop. 3.2.10]{Ma}). Since $\Gm_1$ and $\Gm_2$ are commensurable, we therefore see that 
$F'=F$ (using \cite[VIII, Prop. 3.22]{Ma}) and that $\bo{\wt{G}}\cong\bo{(G')\ii}$ 
(using \cite[I, 0.11]{Ma}). Hence $\bo{P\wt{H}}\cong\bo{R}_{F/\B{Q}}\bo{(G')\ii}$, as claimed.
\end{pf}

\appendix
\section{The Hirzebruch proportionality principle}
The following proposition generalizes \cite[Thm. 2.2.8]{Ku}.
\begin{Prop} \label{P:Hirz}
 Let $X_{\Gm}$ be a projective scheme over a $p$-adic field $L$ such that $\quad X\an_{\Gm}\cong\Gm\bs\prod_{i=1}^r
\Omega^d_{L_i}\otimes_{L_i} L$ for $p$-adic subfields $L_1,...,L_r$ of $L$, a natural number $d\geq 2$,
and a torsion-free cocompact lattice $\Gm\subset\prod_{i=1}^r\PGL{d}(L_i)$ (see \rl{quot}).
Then for any positive integers $i_1,...,i_l$ such that $i_1+...+i_l=r(d-1)$ we have
$c_{i_1}...c_{i_l}(T(X_{\Gm}))=\chi_E(\Gm)\cdot c_{i_1}...c_{i_l}(T(\B{P}^{d-1})^r)$, where by
$c_{i_1}...c_{i_l}$ we denote the corresponding Chern number, and by $\chi_E(\Gm)$ the Euler-Poincar\'e
characteristic of $\Gm$.
\end{Prop}

\begin{pf}
The proof is essentially the same as that of Kurihara, so we will content ourselves with a brief sketch.
In addition we give a proof of a key lemma which Kurihara merely states.

Let $R$ and $\kappa\cong\B{F}_q$ be the ring of integers and the residue field of $L$ respectively. 

\begin{Lem} \label{L:Chern}
(\cite[Lem. 4.6.1]{Ku}) Let $Y$ be a flat projective scheme over $R$ which is locally a complete intersection
over $R$ of relative dimension $n$. Then for any positive integers $i_1,...,i_l$ such that $i_1+...+i_l=n$ we have
\begin{equation} \label{E:Chern}
c_{i_1}...c_{i_l}(T_{Y_{\eta}/L})=\sum_E\text{length}(\C{O}_{Y_0,e})c_{i_1}...c_{i_l}(\imath_E^*(T_{Y_0/\kappa})), 
\end{equation}
where $Y_{\eta}$ (resp. $Y_0$) is the generic (resp. the special) fiber of $Y$, $E$ runs over 
irreducible components of $Y_0$, $e$ denotes the generic point of $E$, $\imath_E$ denotes 
the closed immersion $E\hra Y_0$, and $T_{Y_{\eta}/L}$ (resp. $T_{Y_0}/\kappa$) denotes the virtual tangent bundle of 
$Y_{\eta}/L$ (resp.  of $Y_0/\kappa$).
\end{Lem}
 The virtual tangent bundle of a locally complete intersection morphism $Z_1\to Z_2$ is a certain
 element in the Grothendieck group $K^0(Z_1)$ of vector bundles on $Z_1$ defined for example
in \cite[App. B.7.6]{Fu}.
\begin{pf}
Let $i:\Spec\,\kappa\hra\Spec\,R$ and $j:\Spec\,L\hra\Spec\,R$ be the natural closed regular and open 
embeddings respectively. Let $T_{Y/R}$ be the virtual tangent bundle of $Y/R$. Then $T_{Y_{\eta}/L}=j^*T_{Y/R}$ and 
$T_{{Y_0}/\kappa}=i^*T_{Y/R}$.
Therefore for the proof of the lemma it will suffice to show that for any $l$ classes
$F_1,...,F_l\in K^0(Y)$ we have

$$c_{i_1}(j^*F_1)...c_{i_l}(j^*F_l)=\sum_E\text{length}(\C{O}_{Y_0,e})c_{i_1}(\imath_E^*(i^*F_1))...c_{i_l}
(\imath_E^*(i^*F_l)).$$

Since Grothendieck group is generated by vector bundles we may suppose that all the $F_i$'s are vector bundles. 
We will use the notation and results from \cite{Fu}. 

Recall that the the map $j^*$ from cycles on $Y$ to cycles on $Y_{\eta}$ is surjective and that 
the specialization map  $\sigma$ from cycles on $Y_{\eta}$ to cycles on $Y_0$ defined by $\sigma(j^*\beta):=i^!(\beta)$
is well-defined (see \cite[20.3]{Fu}). Then $\sigma[Y_{\eta}]=[Y_0]=\sum_E\text{length}(\C{O}_{Y_0,e})(\imath_E)_*[E]$ 
and $\deg(\sigma(\al))=\deg(\al)$ for each $0$-cycle $\al$ on $Y_{\eta}$.
Moreover, by the commutativity
of Chern classes with Gysin maps (\cite[Prop. 6.3]{Fu}) and with a flat pull-back (\cite[Thm. 3.2.(d)]{Fu})
we obtain that
$$\sigma(c_m(j^*F)\cap\al)=\sigma(c_m(j^*F)\cap j^*\beta)=\sigma(j^*(c_m(F)\cap\beta))=$$
$$=i^!(c_m(F)\cap\beta)=c_m(i^*F\cap i^!\beta)=c_m(i^*F\cap\sigma(\al))$$ 
for each vector bundle $F$ on $Y$, each cycle $\al=j^*(\beta)$
on $Y_{\eta}$ and each $m\in\B{N}$.
By induction this implies that for any $l$ vector bundles $F_1,...,F_l$ on $Y$ we have 
$$\sigma(c_{i_1}(j^*F_1)\cap...\cap c_{i_l}(j^*F_l)\cap[Y_{\eta}])=\sum_E\text{length}(\C{O}_{Y_0,e})c_{i_1}(i^*F_1)\cap...\cap c_{i_l}(i^*F_l)\cap(\imath_E)_*[E].$$

Using the projection formula \cite[Thm. 3.2(c)]{Fu} and taking degrees of both sides, we obtain the required equality.
\end{pf}

In our case, $X_{\Gm}$ is a generic fiber of a scheme $\wt{X}_{\Gm}$ satisfying the condition of the lemma 
(see \cite[App. 1]{Mus}). Each two irreducible components of the special fiber $\wt{X}_{\Gm,0}$ have isomorphic
neighborhoods, which are moreover  isomorphic to the product of $r$ neighborhoods of components in the case of 
one factor. Therefore all the  summands from the right hand side of (\ref{E:Chern}) are equal.
Using \cite[Thm. 7]{Se} we see that the number of irreducible components in $\wt{X}_{\Gm,0}$ is
$$\Bigl(\frac{d}{(1-q)(1-q^2)\cdot...\cdot(1-q^{d-1})}\Bigr)^r\chi_E(\Gm),$$
hence to prove the proposition it remains to show that
$$c_{i_1}...c_{i_l}(T(\B{P}^{d-1})^r)=\Bigl(\frac{d}{(1-q)(1-q^2)\cdot...\cdot(1-q^{d-1})}\Bigr)^r 
c_{i_1}...c_{i_l}(\imath_E^*(T_{\wt{X}_{\Gm,0}/\B{F}_q})).$$
By the multiplicativity of Chern classes this follows from the the corresponding formula \cite[(4.6.3)]{Ku} 
in the case of one factor.
\end{pf}

\begin{Rem} \label{R:Hirz}
The proposition can be considered as a numerical evidence for the First Main Theorem in two different ways:

 1) It gives evidence (using the classical Hirzebruch proportionality \cite{Hi})
 for \rco{cover} and thus to the results of Step 2 of section 3 (see \rr{equiv}). Moreover in the case $r=1$ it actually
 implies them (see \rr{Yau}).

 2) Using the proposition one can show (by the arguments of \cite[2.6]{V}) that the volumes (=top Chern classes) of
 Shimura variety and of the corresponding $p$-adically uniformized variety are equal.

\end{Rem}


\begin{thebibliography}{99}

\bibitem[BC]{BC}
  J. -F. Boutot et H. Carayol, {\em Uniformisation {$p$}-adique des courbes de Shimura: Les th\'eor\`emes de Cherednik et de Drinfeld}, Ast\'erisque {\bf 196--197} (1991) 45--156.


\bibitem[Be]{Be}
 V. G. Berkovich, {\em The automorphism group of the Drinfeld half-plane}, C.R. Acad. Sci. Paris, Ser. I, Math. 
{\bf 321} (1995) 1127--1132.



\bibitem[Bo]{Bo}
A. Borel, {\em Linear algebraic groups}, Springer, Berlin, 1991.

\bibitem[BZ]{BZ}
  J. -F. Boutot and Th. Zink, {\em The $p$-adic uniformization of Shimura curves}, Preprint des SFB 343, Bielefeld 1995.



\bibitem[CF]{CF}
 J. W. S. Cassels and A. Fr\"olich (eds.), {\em Algebraic number theory}, Academic Press, New York,
1967.


\bibitem[Ch]{Ch}
  I. V. Cherednik, {\em Uniformization of algebraic curves by discrete arithmetic subgroups of {$\PGL{2}(K_w)$} with 
compact quotients}, Math. USSR Sb. {\bf 29} (1976) 55--85.

\bibitem[Cl]{Cl}
 L. Clozel, {\em Repr\'esentations galoisiennes associ\'ees aux repr\'esentations
automorphes autoduales de {$\GL{}(n)$}}, Inst. Hautes \'Etudes Sci. Publ. Math. {\bf 73} (1991) 77--145.

\bibitem[De1]{De}
 P. Deligne, {\em Travaux de Shimura}, Sem. Bourbaki {\bf 389}, Lecture Notes in
Math. {\bf 244}, Springer, Berlin, 1971, 123--165.

\bibitem[De2]{De2}
  \bysame, {\em Vari\'etes de Shimura: Interpr\'etation modulaire, et techniques de 
construction de mod\'eles canoniques}, Proc. Sympos. Pure Math. {\bf 33}, Part II, 
Amer. Math. Soc. (1979) 247--290.



\bibitem[Dr1]{Dr1}
  V. G. Drinfeld, {\em Elliptic modules}, Math. USSR Sb. {\bf 23} (1974) 561--592.

\bibitem[Dr2]{Dr}
 \bysame, {\em Coverings of {$p$}-adic symmetric regions}, Functional Anal. Appl. {\bf 10} (1976) 107--115.

\bibitem[Fu]{Fu}
  W. Fulton, {\em Intersection Theory}, Springer-Verlag, Berlin, 1984.


\bibitem[Gu]{Gu}
 R. Gunning, {\em Introduction to holomorphic functions of several variables},
Wadsworth, Brooks/Cole, Math. Series, Belmont, Calif., 1990.




\bibitem[He]{He}
 S. Helgason, {\em Differential geometry and symmetric spaces}, Academic Press, New York,
 1962.

\bibitem[Hi]{Hi}
 F. Hirzebruch, {\em Automorphe Formen und der Satz von Riemann-Roch},
Sympos. Internat. Topol. Alg., Univ. de Mexico, 1958, 129--144.


\bibitem[Kn]{Kn}
 A. Knapp, {\em Lie groups beyond the introduction}, Birkh\"auser, Basel, 1996.


\bibitem[Kob]{Kob}
 S. Kobayashi, {\em Transformation groups in differential geometry}, Springer, 
 Berlin, 1972.
 
\bibitem[Kot]{Ko}
 R. E. Kottwitz, {\em Tamagawa numbers}, Ann. of Math. {\bf 127} (1988) 629--646.

\bibitem[Ku]{Ku}
  A. Kurihara, {\em Construction of {$p$}-adic unit balls and the Hirzebruch 
proportionality}, Amer. J. Math. {\bf 102} (1980) 565--648.

\bibitem[La]{La}
 S. Lang, {\em Introduction to complex hyperbolic spaces}, Springer, Berlin, 1987.

\bibitem[Ma]{Ma}
 G. A. Margulis, {\em Discrete subgroups of semisimple Lie groups}, 
Springer, Berlin, 1991.

\bibitem[Mi1]{Mi1}
 J. Milne, {\em Canonical models of (mixed) Shimura varieties}, in 
{\it Automorphic Forms, Shimura Varieties and {$L$}-functions}, Vol {\bf 1}, 
(L. Clozel and J. Milne eds.), Academic Press, New York, 1990, 283--414.


\bibitem[Mi2]{Mi2}
 \bysame, {\em Kazhdan's theorem on arithmetic varieties}, handwritten notes, 1984.

\bibitem[Mi3]{Mi3}
 \bysame, {\em The points on a Shimura variety modulo a prime of good reduction}, 
in {\em Zeta Functions of Picard Modular Surfaces}, (R. Langlands and D. Ramakrisnan eds.), Univ. Montreal, 1992, 151--253.

\bibitem[Mi4]{Mi4}
 \bysame, {\em Shimura varieties and motives}, Proc. Sympos. Pure Math. {\bf 55}, Part II, Amer. Math. Soc. (1994) 447--523.

\bibitem[Mum]{Mum}
 D. Mumford, {\em An analytic construction of degenerating curves over complete
 local rings}, Compositio Math. {\bf 24} (1972) 129--174.


\bibitem[Mus]{Mus}
 G. A. Mustafin, {\em Nonarchimedean uniformization}, Math. USSR Sb. {\bf 34} 
 (1978) 187--214.

\bibitem[NR]{NR}
 M. Nori and M. Raghunathan, {\em On conjugation of locally symmetric arithmetic 
 varieties}, Proc. Indo-French Conference on Geometry, 
(S. Ramanan and A. Beauville eds.), Hindustan Book Agency, 1993, 111--122.


\bibitem[PR]{PR}
 V. Platonov and A. Rapinchuk, {\em Algebraic groups and number theory},
Academic Press, New York, 1994.

\bibitem[PV]{PV}
 M. van der Put and H. Voskuil, {\em Symmetric spaces associated to split algebraic 
groups over a local field}, J. Reine Angew. Math. {\bf 433} (1992) 69--100.


\bibitem[Ra]{Ra}
  M. Rapoport, {\em On the bad reduction of Shimura varieties}, in {\em Automorphic 
Forms, Shimura Varieties and {$L$}-Functions}, Vol {\bf 2}, (L. Clozel and J. Milne
eds.), Academic Press, New York, 1990, 253--321.

\bibitem[RZ1]{RZ1}
 \bysame and Th. Zink, {\em \"{U}ber die okale Zetafunktion von Shimuraariet\"aten}, 
Invent. Math. {\bf 68} (1982) 21--101. 

\bibitem[RZ2]{RZ2}
 \bysame, {\em Period spaces for $p$-divisible groups}, Annals of Math. Stud. {\bf 141}, 
Princeton University Press, 1996.



\bibitem[Sha]{Sh}
 I. R. Shafarevich, {\em Basic algebraic geometry}, Springer, Berlin, 1974.

\bibitem[Shi]{Shi}
 G. Shimura, {\em Introduction to the arithmetic theory of automorphic functions}, 
Princeton Univ. Press, Princeton, 1971.
   
\bibitem[Se]{Se}
  J. -P. Serre, {\em Cohomologie des Groupes Discrets},  Prospects in 
Mathematics, Ann. of Math. Stud. {\bf 70}, Princeton Univ. Press, (1971) 77--169.





 


\bibitem[Ti]{Ti}
  J. Tits, {\em Classification of algebraic semisimple groups}, Proc. Sympos. Pure 
Math. {\bf 9}, Amer. Math. Soc. (1966) 32--62.

\bibitem[V]{V}
   Y. Varshavsky, {\em $P$-adic uniformization of unitary Shimura
   varieties}, Inst. Hautes \'Etudes Sci. Publ. Math. {\bf 87} 
(1998) 57--119.




\bibitem[Ya1]{Ya1}
 S. -T. Yau, {\em Calabi's conjecture and some results in algebraic geometry},
Proc. Nat. Acad. Sci. U.S.A. {\bf 74}, (1977), pp. 1798--1799.

\bibitem[Ya2]{Ya}
 S. -T. Yau, {\em On the Ricci curvature of a compact K\"{a}hler manifold and the
complex Monge-Amp\`ere equation, I}, Comm. Pure Appl. Math. {\bf 31}
(1978) 339--411.

\end{thebibliography}
\end{document}